\documentclass[11pt]{pjm}
\mathchardef\varSigma="0106
\usepackage{amssymb,latexsym}
\begin{document} 
\font\tensans=cmss10
\def\hs{\hskip.7pt}
\def\nh{\hskip-.7pt}
\def\lj{{\bf[}}
\def\rj{{\bf]}}
\def\bbR{{\bold R}}
\def\bbC{{\bold C}}
\def\bbH{{\bold H}}
\def\bbK{{\bold K}}
\def\bbF{{\bold F}}
\def\cp{{\bold C}{\rm P}\hskip.4pt}
\def\hp{{\bold H}{\rm P}}
\def\kp{{\bold K}{\rm P}}
\def\proj{{\rm pr}\hs}
\def\emp{{\tensans\O}}
\def\bbZ{{\bold Z}}
\def\h{{\rm Hom}\hskip.5pt}
\def\hr{{\rm Hom}_{\hskip.4pt\bbR}}
\def\hc{{\rm Hom}_{\hskip.4pt\bbC}}
\def\hk{{\rm Hom}_{\hskip.4pt\bbK}}
\def\hh{{\rm Hom}_{\hskip.4pt\bbH}}
\def\sug{{\rm SU}\hskip.5pt}
\def\ug{{\rm U}\hskip.5pt}
\def\gy{G_2^+\hskip-1.8pt\y}
\def\gv{G_2^+(\mathcal{V})}
\def\tgy{T\hskip1pt[G_2^+\hskip-1.6pt\y]}
\def\vrad{v^{\hs\rm rad}}
\def\vtng{v^{\hs\rm tng}}
\def\tr{totally real}
\def\tri{totally real immersion}
\def\psh{pseu\-do\-hol\-o\-mor\-phic}
\def\mfd{-man\-i\-fold}
\def\sbms{submanifolds}
\def\ac{almost complex}
\def\acsu{almost complex surface}
\def\acst{almost complex structure}
\def\y{{M}}
\def\x{\varSigma}
\def\vg{\varGamma}
\def\vt{\varTheta}
\def\dvt{d\hskip.4pt\vt}
\def\zr{Z}
\def\ve{\varepsilon}
\def\pep{Q\hs}
\def\pem{\hskip1.2pt\overline{\hskip-1.2ptQ}\hs}
\def\pw{P(W)}
\def\pv{P(\mathcal{V})}
\def\pmw{\hskip1.2pt\overline{\hskip-1.2ptP}(W)}
\def\pmv{\hskip1.2pt\overline{\hskip-1.2ptP}(\mathcal{V})}
\def\hf{F}
\def\fo{\varOmega}
\def\ap{\alpha}
\def\be{\beta}
\def\ga{\gamma}
\def\spanc{{\rm Span}_{\hskip.5pt\bbC\hskip.3pt}}
\def\spanr{{\rm Span}_{\hskip.5pt\bbR\hskip.3pt}}
\def\dimr{\dim_{\hskip.4pt\bbR\hskip-1.2pt}}
\def\dimc{\dim_{\hskip.4pt\bbC\hskip-1.2pt}}
\def\dimh{\dim_{\hskip.4pt\bbH\hskip-1.2pt}}
\def\dimk{\dim_{\hskip.4pt\bbK\hskip-1.2pt}}
\def\detr{{\rm det}_\bbR}
\def\detc{{\rm det}_\bbC}
\def\tracer{{\rm Trace}_{\hskip.4pt\bbR}}
\def\df{d\hskip-.8ptf}
\def\Spinc{{\rm Spin}$^{\hskip.3pt{\rm c}}$}
\def\bold{\bf}
\def\ce{$^{\hskip.5pt{\rm c}}$}
\def\cf{$^{\hskip.5pt{\rm c}}(4)$}
\def\spingeo{spin\cf-ge\-om\-e\-try}
\def\spst{{\rm spin}$^{\hskip.3pt{\rm c}}$-struc\-ture}
\def\csg{\mathcal{S}}
\def\sg{\sigma}
\def\cka{\mathcal{K}}
\def\ka{\kappa}
\def\la{\lambda}
\def\diml{-di\-men\-sion\-al}
\title[Spin\ce-man\-i\-folds with Higgs fields]{Immersions of surfaces in 
spin$^{\hskip-.2pt{\rm c}\hskip.8pt}$-man\-i\-folds with Higgs fields} 
  
  
\author[A. Derdzinski]{Andrzej Derdzinski} 
\address{Dept. of Mathematics\\
The Ohio State University\\
Columbus, OH 43210, USA} 
\email{andrzej@math.ohio-state.edu} 
\urladdr{http://www.math.ohio-state.edu/\~{}andrzej}
\author[T. Januszkiewicz]{Tadeusz Januszkiewicz}
\address{Mathematical Institute\\
Wroc\l aw University\\
pl. Grunwaldzki 2/4\\ 
50-384 Wroc\l aw, Poland\\
and\\
Mathematical Institute, Polish Academy of Sciences} 
\urladdr{http://www.math.uni.wroc.pl/\~{}tjan} 

\thanks{The second author was partially supported by the Ohio State University 
Research Foundation (OSURF) and Poland's Committee of Scientific Research 
(KBN), under KBN grant no.~2~P03A~035~20}

\newtheorem{thm}{Theorem}[section] 
\newtheorem{prop}[thm]{Proposition} 
\newtheorem{lem}[thm]{Lemma} 
\newtheorem{cor}[thm]{Corollary} 
  
\theoremstyle{definition} 
  
\newtheorem{defn}[thm]{Definition} 
\newtheorem{notation}[thm]{Notation} 
\newtheorem{example}[thm]{Example} 
\newtheorem{conj}[thm]{Conjecture} 
\newtheorem{prob}[thm]{Problem} 
  
\theoremstyle{remark} 
  
\newtheorem{rem}[thm]{Remark} 
  
\date{Received date / Revised version date} 

\begin{abstract}
The usual assumption one makes to define totally real or 
pseu\-do\-hol\-o\-mor\-phic immersions of real surfaces in a four-man\-i\-fold 
$\,\y\,$ is that $\,\y\,$ carries a fixed almost complex structure. We extend 
both definitions to a more general case of a spin\ce-man\-i\-fold $\,\y\,$ 
with a `Higgs field', that is, a generic smooth section of the positive 
half-spinor bundle, and describe all pseu\-do\-hol\-o\-mor\-phic immersions of 
closed surfaces in the four-di\-men\-sion\-al sphere endowed with a standard 
Higgs field.
\end{abstract} 
  
\maketitle 

\section{Introduction}\label{intro} 
Almost complex structures on real manifolds of dimension $\,2n\,$ are 
well-known to be, essentially, a special case of \spst s. This amounts to a 
specific Lie-group embedding 
$\,\hs\ug(n)\to$\hskip4ptSpin$^{\hskip.5pt{\rm c}}(2n)\,$ 
(see~\cite[p.~392]{lawson-michelsohn} and Remark~\ref{reemb} below). The 
present paper deals with the case $\,n=2$. The relation just mentioned then 
can also be couched in the homotopy theorists' language: {\em for a compact 
four\mfd\/ $\,\y\,$ with a fixed\/ $\,{\rm CW}$-decomposition, a \spst\ over\/ 
$\,\y\,$ is nothing else than an \acst\ on the\/ $\,2$-skeleton of\/ 
$\,\y\nh$, admitting an extension to its\/ $\,3$-skeleton\/{\rm;} the 
extension itself is not a part of the data.} (Kirby \cite{kirby} attributes 
the italicized comment to Brown.)

Our approach is explicitly geometric and proceeds as follows. Given an \acst\ 
$\,J\,$ on a $\,4$\mfd\ $\,\y$, one can always choose a Riemannian metric 
$\,g\,$ compatible with $\,J$, thus replacing $\,J\,$ by an almost Hermitian 
structure $\,(J,g)\,$ on $\,\y$. The latter may in turn be treated as a \spst\ 
on $\,\y\,$ with a fixed unit $\,C^\infty$ section $\,\psi\,$ of its positive 
half-spinor bundle $\,\sg^+$ (cf.\ \S\ref{acsu}). Complex points of an 
immersion in $\,\y\,$ of any oriented real surface $\,\x\,$ then can be 
directly described in terms of $\,\psi\,$ rather than $\,J$. Such a 
description is presented in the initial part of the paper, culminating in 
Theorem~\ref{thuno}. Moreover, when $\,\x\,$ is compact, the net total number 
of complex points for a generic immersion $\,\x\to\y\,$ depends just on the 
underlying \spst, not even on $\,\psi\,$ (see Corollary~\ref{corol}).

Imitating the above description, we may define complex points of a 
real-surface immersion $\,\x\to\y$, where $\,\y\,$ is a $\,4$\mfd\ with a 
\spst\ and {\em any\/} fixed section $\,\psi\,$ of $\,\sg^+$. This raises the 
question of finding a suitable additional requirement on $\,\psi$, for which 
the corresponding notion of a complex point and the resulting classes of \tr\ 
and \psh\ immersions would be less restrictive than for \acst s, yet still 
interesting and reasonable. (As an ``unreasonable'' example, choose $\,\psi\,$ 
to be the zero section; then {\em every\/} immersion is \psh.)

The definition we propose, in \S\ref{hggs}, is that of a {\em Higgs field\/} 
for a given \spst\ on a $\,4$\mfd, which is a $\,C^\infty$ section $\,\psi\,$ 
of $\,\sg^+$ transverse to the zero section, and defined only up to a positive 
functional factor.

The notions of {\em \tr\/} and {\em \psh\/} immersions of oriented real 
surfaces $\,\x\,$ now have straightforward extensions to the case where the 
receiving $\,4$\mfd\ $\,\y\,$ is endowed with a \spst\ and a Higgs field. We 
observe, in \S\ref{trei}, that the former immersions have properties analogous 
to those in the almost-complex case, and classify the latter ones when 
$\,\y\,$ is $\,S^{\hs4}$ with a ``standard'' Higgs field and $\,\x\,$ is a 
closed surface (see \S\ref{rsps}).

Unlike \acst s, Higgs fields exist in every \spst, which in turn exists on 
every $\,4$\mfd. On the other hand, pseu\-do\-hol\-o\-mor\-phic immersions of 
surfaces have been used, with great success, to probe the topology of 
symplectic manifolds, especially in dimension $\,4$. One may therefore wonder 
if they could similarly be employed as a tool for studying more general 
$\,4$\mfd s with Higgs fields.

\section{Preliminaries}\label{prel}
'Planes' and 'lines' always mean {\em vector\/} spaces. A real/com\-plex 
vector space $\,V\hs$ is called {\em Euclidean/Hermitian\/} if 
$\,\dim V<\infty\,$ and $\,V\hs$ carries a fixed pos\-i\-tive-def\-i\-nite 
inner product (always denoted $\,\langle\,,\rangle$, with the symbol 
$\,|\,\,|\,$ used for its associated norm). Such $\,\langle\,,\rangle\,$ is 
uniquely determined by $\,\,{\rm Re}\,\langle\,,\rangle\,$ or $\,|\,\,|$.
\begin{rem}\label{relin}We identify oriented (real) Euclidean planes with 
(complex) Hermitian lines, so that both inner products have the same real part 
and norm, and multiplication by $\,i\,$ is the positive rotation by the angle 
$\,\pi/2$.
\end{rem}
All manifolds are assumed connected except when stated otherwise, while most 
real manifolds we deal with are oriented and of class $\,C^\infty$. The 
orientations for orthogonal complements of real vector subspaces Cartesian 
products, total spaces of locally trivial bundles, preimages of regular values 
of mappings, and (cf.\ (a) in \S\ref{trsv}) zero sets of transverse sections 
in vector bundles, are all obtained using the direct-sum convention; since the 
real dimensions involved are all even, no ambiguity arises even if the order 
of the summands is not specified. This is consistent with the convention which 
treats (almost) complex manifolds as {\em oriented\/} real manifolds, 
declaring that
\begin{equation}
\aligned
&(e_1,ie_1,\dots,e_n,ie_n)\hskip5.5pt{\rm is\ a\ positive{\textrm-}oriented\ 
real\ basis\ of\ a}\\
&{\rm complex\ vector\ space\ with\ a\ complex\ basis}\hskip6.5pt
(e_1,\dots,e_n){\rm.}\endaligned\label{ori}
\end{equation}
We also deal with {\em inner-product spaces}, including {\em normed lines}, 
over the field $\,\bbH\,$ of quaternions. The $\,\bbH$-sesqui\-linearity 
requirement imposed on an $\,\bbH$-val\-ued inner product 
$\,\langle\,,\rangle\,$ in a quaternion vector space $\,W\,$ with 
$\,\dim W<\infty\,$ includes the condition 
$\,\langle px,qy\rangle=p\hs\langle x,y\rangle\hs\overline q\,$ for 
$\,x,y\in W\,$ and $\,p,q\in\bbH\hs$. Using the inclusion  
$\,\bbC=\,\spanr(1,\hs{\bf i})\,\subset\,\bbH\hs$, we will treat any such 
$\,W\,$ as a complex Hermitian space, with the $\,\bbC$-val\-ued inner product 
having the same real part (or, equivalently, the same associated norm) as the 
original $\,\langle\,,\rangle$.

Denoting $\,\bbK\,$ any of the scalar fields $\,\bbR\hs,\bbC\hs,\bbH\hs$, 
we let $\,G_m(W)\,$ stand for the Grassmannian manifold of all subspaces 
$\,L\,$ with $\,\dimk L=m\,$ in a given vector space $\,W\,$ over $\,\bbK\hs$, 
where $\,1\le m\le\dim W<\infty$. Each tangent space $\,T_L[G_m(W)]\,$ of 
$\,G_m(W)\,$ then has a canonical real-i\-so\-mor\-phic identification 
$\,T_L[G_m(W)]=\,\hk(L,W/L)$. Namely, $\,A\in\,\hk(L,W/L)\,$ corresponds to 
$\,d\pi_{\bf w}[RA\hs{\bf w}]\in T_L[G_m(W)]\,$ for any basis 
$\,{\bf w}=(w_1,\dots,w_m)\,$ of $\,L\,$ and any right inverse $\,R\,$ of the 
projection $\,W\to W/L$. Here $\,\pi\,$ is the standard projection onto 
$\,G_m(W)\,$ of the {\em Stiefel manifold\/} $\,\hs{\rm St}_{\hs m}(W)$, i.e., 
the open set in the $\,m$th Cartesian power $\,W^m$ formed by all linearly 
independent systems $\,{\bf w}=(w_1,\dots,w_m)$, and 
$\,RA\hs{\bf w}=(RAw_1,\dots,RAw_m)$, so that 
$\,RA\hs{\bf w}\in W^m=T_{\bf w}[W^m]=T_{\bf w}[{\rm St}_m(W)]$.
\begin{rem}\label{regra}Thus, writing any given vector in $\,T_L[G_m(W)]\,$ as 
$\,d\pi_{\bf w}{\bf v}\,$ with a fixed 
$\,{\bf w}\in\pi^{-1}(L)\subset\hs{\rm St}_{\hs m}(W)\,$ and a suitable 
$\,{\bf v}=(v_1,\dots,v_m)\in W^m=T_{\bf w}[{\rm St}_m(W)]$, we can define 
$\,A\in\,\hk(L,W/L)\,$ associated as above with 
$\,d\pi_{\bf w}{\bf v}\in T_L[G_m(W)]\,$ to be the operator sending the basis 
$\,(w_1,\dots,w_m)\,$ of $\,L\,$ onto the system $\,(v_1+L,\dots,v_m+L)\,$ in 
$\,W/L$.
\end{rem}
Next, let $\,\pw\approx\hs\kp^{n-1}$ be the projective space of all lines in a 
vector space $\,W\,$ over $\,\bbK\,$ with $\,1\le\dim W<\infty$. As 
$\,\pw=G_1(W)$, the above identification $\,T_L[G_m(W)]=\,\hk(L,W/L)\,$ now 
becomes
\begin{equation}
T_L[\pw]\,=\,\hs\hk(L,W/L)\qquad{\rm for\ every}\hskip8ptL\in\pw\hs.\label{tlp}
\end{equation}
\begin{rem}\label{redvt}Given $\,W,\bbK\,$ as above, a codimension-one 
subspace $\,V\hs$ of $\,W$, and a vector $\,u\in W\smallsetminus V$, let 
$\,\y'$ be the complement of $\,P(V)\,$ in $\,\pw$. The mapping 
$\,\vt:V\to\y'$ defined by $\,\vt(y)=\bbK w\,$ with $\,w=y+u\,$ then clearly 
is a $\,C^\infty$ diffeomorphism, while $\,\vt^{-1}$ represents a standard 
{\em projective coordinate system\/} in $\,\pw$. Also, let $\,A=\dvt_yv\,$ for 
any $\,y\in V$ and $\,v\in V=T_yV$. Using (\ref{tlp}) with $\,L=\vt(y)\,$ to 
treat $\,A\,$ as an operator $\,L\to W/L$, we then have $\,Aw=v+L$, where 
$\,w=y+u$.

This is clear from Remark~\ref{regra} for $\,m=1$, as $\,\vt\,$ is the 
restriction to $\,V$ of the translation by $\,u\,$ followed by the projection 
$\,\pi:W\smallsetminus\{0\}\to \pw$, and so $\,\dvt_yv=d\pi_wv\,$ by the chain 
rule.
\end{rem}

\section{The Grassmannian of oriented planes}\label{grop}
Let $\,G_2^+(W)\,$ be the Grassmannian of real oriented planes in a real 
vector space $\,W\hs$ with $\,\dim W<\infty$, so that $\,G_2^+(W)\,$ is a 
two-fold covering manifold of $\,G_2(W)\hs$ defined as in \S\ref{prel} for 
$\,\bbK=\bbR\hs$. Any inner product in $\,W\,$ makes $\,G_2^+(W)\,$ an \ac\ 
manifold, and so, by (\ref{ori}), $\,G_2^+(W)\,$ {\em is naturally oriented}. 
In fact, an oriented Euclidean plane $\,\mathcal{T}\in G_2^+(W)\,$ is a 
complex line (Remark~\ref{relin}), which turns the tangent space of 
$\,G_2^+(W)\,$ at $\,\mathcal{T}$, i.e., 
$\,\hs\hr(\mathcal{T},W/\mathcal{T})\,$ (see \S\ref{prel}), into a complex 
vector space.
\begin{rem}\label{reppg}If, in addition, $\,W\,$ itself is the underlying real 
space of a complex vector space, $\,G_2^+(W)\hs$ contains two disjoint, 
embedded complex manifolds $\,\pw,\hs\pmw$, which are copies of 
the complex projective space of $\,W\,$ consisting of all complex lines in 
$\,W\,$ treated as real planes oriented, in the case of $\,\pw$, by the 
complex-line orientation with (\ref{ori}), or, for $\,\pmw$, by its opposite. 
The \acst\ in $\,G_2^+(W)$, obtained as above from a Euclidean inner product 
in $\,W\,$ which is the real part of a Hermitian inner product, then makes the 
differential of the inclusion $\,\pw\to G_2^+(W)\,$ (or, $\,\pmw\to G_2^+(W)$) 
complex-linear (or, respectively, antilinear) at each point. This is clear 
from (\ref{tlp}) and its analogue for $\,T_L[G_m(W)]\,$ in \S\ref{prel}.
\end{rem}
For a $\,4$\mfd\ $\,\y$, the $\,8$\diml\ Grassmannian manifold $\,\gy\,$ is 
the total space of the bundle over $\,\y\,$ with the fibres $\,G_2^+(T_y\y)$, 
$\,y\in\y$, so that, as a set, 
$\,\gy=\{(y,\mathcal{T}):y\in\y\hskip3pt{\rm and}
\hskip3pt\mathcal{T}\in G_2^+(T_y\y)\}$. The {\em Gauss mapping\/} 
$\,\hf:\x\to\gy\,$ of any immersion 
$\,f:\x\to\y\,$ of an oriented real surface $\,\x\,$ then is given by 
$\,\hf(x)=(f(x),\df_x(T_x\x))$. If $\,M\,$ is oriented, so is $\,\gy$, due to 
the natural orientations of its fibres (see above).

If $\,\y\,$ is the underlying oriented $\,4$\mfd\ of an \acsu, cf.\ 
(\ref{ori}), the oriented $\,8$\diml\ manifold $\,\gy\,$ has two distinguished 
$\,6$\diml\ \sbms\ $\,\pep\,$ and $\,\pem$, which form total spaces of 
$\,\cp^1$-bundles over $\,\y$. Their fibres over any $\,y\in\y\,$ are 
$\,\pw\,$ and $\,\pmw$, defined as in Remark~\ref{reppg} for $\,W=T_y\y$. Both 
$\,\pep\,$ and $\,\pem\,$ are orientable. We will always treat them as 
oriented manifold, choosing, however, the {\em opposites\/} of their natural 
orientations. Specifically, the orientations we use are the direct sums of the 
orientation with (\ref{ori}) on the base $\,\y\,$ and the orientations of the 
fibres $\,\pw\,$ and $\,\pmw$, with $\,W=T_y\y$, which are opposite to their 
natural orientations of complex projective lines.

\section{Complex points of immersed real surfaces}\label{cpts} 
Any real plane $\,\mathcal{T}\,$ in a complex plane $\,V\hs$ (see 
\S\ref{prel}) either is a complex line, or is {\em totally real\/} in the 
sense that $\,\,\spanc\mathcal{T}=V$. 
 
Let $\,\y\,$ be an \acsu\ ($\dimr\y=4$). An immersion $\,f:\x\to\y\,$ of a 
real surface $\,\x\,$ is said to have a {\em complex point\/} at $\,x\in\x\,$ 
if $\,\tau_x=\df_x(T_x\x)\,$ is a complex line in $\,T_{f(x)}\y$. One calls 
$\,f\,$ {\em totally real\/} or {\em \psh\/} if it has no complex points or, 
respectively, only complex points. Thus, $\,f\,$ is totally real if and only 
if $\,\tau_x$ is totally real in $\,T_{f(x)}\y\,$ for every $\,x\in\x$. On the 
other hand, if $\,f\,$ is \psh, $\,\x\,$ acquires a natural orientation, 
pulled back to $\,T_x\x\,$ from $\,\tau_x$ by $\,\df_x$ for every $\,x\in\x$. 
We say that $\,f\,$ is a \psh\ immersion in $\,\y\,$ of an {\em oriented\/} 
real surface $\,\x\,$ if this orientation coincides with the one prescribed in 
$\,\x\,$ or, equivalently, if $\,\hf(\x)\subset\pep\,$ (with $\,\hf,\pep\,$ as 
at the end of \S\ref{grop}). See \cite{chen-ogiue}, \cite{gromov}.

The {\em determinant bundle\/} $\,\detr\hs\eta\,$ (or, $\,\,\detc\hs\eta$) of 
a rank $\,k\,$ real/com\-plex vector bundle $\,\eta\,$ is its highest exterior 
power $\,\eta^{\hs\wedge k}$. For an immersion $\,f\,$ of a real surface 
$\,\x\,$ in an \acsu\ $\,\y$,
\begin{equation}
f^*c_1(\ka)\,=\,0\hskip7pt{\rm in}\hskip6ptH^2(\x,\bbZ)
\hskip8pt{\rm if}\hskip7ptf\hskip4pt{\rm is\ totally\ real\ and}
\hskip5pt\ka=\detc\hs T\y\hs.\label{tre}
\end{equation}
One then also has a natural o\-ri\-en\-ta\-tion-re\-vers\-ing vector-bundle 
isomorphism between the tangent bundle $\,T\x\,$ and the normal bundle 
$\,\nu_f$, that is,
\begin{equation}
\nu_f\,=\,\,\overline{T\x}\hskip15pt{\rm whenever}\hskip10ptf\hskip7pt
{\rm is\ totally\ real\ and}\hskip7pt\x\hskip7pt{\rm is\ oriented.}\label{tno} 
\end{equation}
These well-known facts follow since 
$\,f^*T\y\,=\,\,\spanc\tau\,=\,\tau\oplus\,i\tau$, where $\,\tau$ stands for 
the subbundle $\,\df(T\x)\,$ of $\,f^*T\y$. Thus, 
$\,\nu_f=i\tau$, which gives (\ref{tno}) as the multiplication by $\,i\,$ is 
o\-ri\-en\-ta\-tion-re\-vers\-ing (by (\ref{ori})), while $\,f^*T\y\,$ 
coincides with the complexification $\,\tau^\bbC$ of $\,\tau\approx T\x$, and 
so $\,f^*[\detc\hs T\y]\,=\,[\detr\hs T\x]^\bbC$. Now (\ref{tre}) is obvious: 
namely, $\,\,\detr\hs T\x\,$ is trivial if $\,\x\,$ is closed and orientable, 
and $\,H^2(\x,\bbZ)=\{0\}\,$ otherwise. 

If $\,\y\,$ is an \acsu, $\,f:\x\to\y\,$ is an {\em arbitrary\/} immersion of 
a {\it closed, oriented\/} real surface, and $\,\,\cdot\,\,$ denotes the 
intersection form in $\,H_2(\gy,\hs\bbZ)$, we may replace (\ref{tre}) by 
the following formula, proved later in Corollary~\ref{corol}: with 
$\,\pep\nh,\pem\,$ defined at the end of \S\ref{grop},
\begin{equation}
\aligned
&\int_\x f^*c_1(\ka)\,=\,([\pep]+[\pem])\cdot\hf_*[\x]\hskip17pt
{\rm for}\hskip11pt\ka\hs=\,\detc\hs T\y\nh,\\
&\hskip2pt{\rm where}\hskip6pt\hf:\x\to\gy\hskip6pt{\rm is\ the\ 
Gauss\ mapping\ of}\hskip5ptf\hskip5pt{\rm(see\ \S\ref{grop}).}\endaligned
\label{gss}
\end{equation}
When $\,f\,$ is \tr, $\,\pep\cap\hf(\x)=\pem\cap\hf(\x)=$\hskip4.5pt\emp\hs, 
and (\ref{gss}) becomes (\ref{tre}).

Expression $\,([\pep]+[\pem])\cdot\hf_*[\x]\,$ in (\ref{gss}) represents the 
{\em net total number of complex points\/} of $\,f$. 
In fact, given an immersion $\,f\,$ of an oriented real surface $\,\x\,$ 
(closed or not) in an \acsu\ $\,\y$, every complex point $\,x\,$ of $\,f\,$ 
is a $\,\pep${\it\hskip-1.4pt-com\-plex point\/} or a 
$\,\pem${\it\hskip-1.4pt-com\-plex point}, in the sense that 
$\,\hf(x)\in\pep\,$ or, respectively, $\,\hf(x)\in\pem$. The {\em index\/} of 
such an immersion $\,f$, at any isolated complex point $\,x\in\x$, is defined 
to be the intersection index that the immersion $\,\hf:\x\to\gy\,$ has, at the 
isolated intersection point $\,x$, with the disconnected oriented $\,6$\mfd\ 
$\,\pep\cup\pem\subset\gy$.

Thus, if $\,\x\,$ is closed and the immersion $\,f\,$ is {\em generic\/} 
(i.e., has only finitely many complex points), 
$\,([\pep]+[\pem])\cdot\hf_*[\x]\,$ 
is the sum of indices of all complex points of $\,f$, while the separate 
contributions corresponding to $\,\pep$\hskip-1.2pt-com\-plex and 
$\,\pem$\hskip-1.2pt-com\-plex points are $\,[\pep]\cdot\hf_*[\x]\,$ and 
$\,[\pem]\cdot\hf_*[\x]$. 

\section{Spin$^{\hskip.3pt{\rm c}}(4)$-geometries}\label{spge}
The terms 'Hermitian plane/line' refer, as in \S\ref{prel}, to complex vector 
spaces.

We treat the quaternion algebra $\,\bbH\,$ as a Hermitian plane, with the 
multiplication by complex scalars declared to be the {\em right\/} quaternion 
multiplication by elements of 
$\,\bbC=\,\spanr(1,\hs{\bf i})\,\subset\,\bbH\hs$, and with the inner product 
$\,\langle\,,\rangle\,$ making the $\,\bbC$-basis $\,1,\,{\bf j}\,\,$ 
orthonormal, so that it corresponds to the norm $\,|\,\,|\,$ with 
$\,|p|^2=p\overline p$. Somewhat surprisingly,
\begin{equation}
1,\,{\bf i}\hs,\,{\bf j}\hs,\,{\bf k}\hskip6pt{\rm form\ a\ }{\it 
negative{\textit-}oriented\/\ }{\rm real\ basis\ of}\hskip6pt\bbH\hs{\rm,}
\label{neg}
\end{equation}
as the $\,\bbC$-basis $\,1,\,{\bf j}\hs\,$ leads to the 
pos\-i\-tive-o\-ri\-ented $\,\bbR$-basis 
$\,1,\,{\bf i}\hs,\,{\bf j}\hs,\,{\bf j}{\bf i}\,$ (see 
(\ref{ori})) and $\,\hs{\bf j}{\bf i}\hs=\hs-\hs{\bf k}\hs$. For 
$\,a,b,x,y,u,v\in\bbC\,$ we have $\,a{\bf j}={\bf j}\hs\overline a$, and so, 
in $\,\bbH\nh$,
\begin{equation}
(a+\hskip.0pt{\bf j}b)(x+\hskip.0pt{\bf j}y)\,
=\,u+\hskip.0pt{\bf j}v\hskip13pt{\rm if\ and\ only\ if}\hskip11pt
\left[\begin{matrix}a&-\,\overline b\,\cr
b&\,\,\overline a\end{matrix}\right]
\left[\begin{matrix}x\cr
y\end{matrix}\right]
=\left[\begin{matrix}u\cr
v\end{matrix}\right]\nh.
\label{mlt}
\end{equation}
\begin{defn}\label{defin}By a spin\cf{\it-ge\-om\-e\-try\/} we mean a triple 
$\,(\csg^+\!,\csg^-\!,\cka)\,$ consisting of two Hermitian planes $\,\csg^\pm$ 
and a Hermitian line $\,\cka$, endowed with two fixed norm-pre\-serv\-ing 
isomorphic identifications $\,[\csg^\pm]^{\wedge2}=\cka$. 

The inner products in $\,[\csg^\pm]^{\wedge2}$, used here, are induced by 
those of $\,\csg^\pm$, via the formula 
$\,|\phi\wedge\chi\hs|^2=|\phi|^2|\chi\hs|^2-|\langle\phi,\chi\rangle|^2$. 
Thus, instead of assuming that such identifications are given, we could 
require that there be skew-sym\-met\-ric bilinear multiplications 
$\,\csg^\pm\times\,\csg^\pm\to\cka\,$ which, written as 
$\,(\phi,\chi)\mapsto\phi\wedge\chi$, satisfy the last formula or, 
equivalently, the condition $\,|\phi\wedge\chi\hs|=|\phi|\,|\chi\hs|\,$ 
whenever $\,\langle\phi,\chi\rangle=0$.
\end{defn}
An example of a \spingeo\ is $\,\csg^\pm=\bbC^2$, $\,\cka=\bbC\,$ with both 
skew-sym\-met\-ric multiplications given by $\,\phi\chi=ad-bc\,$ for 
$\,\phi=(a,b)$, $\,\chi=(c,d)$. Any \spingeo\ $\,(\csg^+\!,\csg^-\!,\cka)\,$ 
is equivalent to this one under an isomorphism obtained by choosing bases 
$\,(\phi^\pm\nh,\chi^\pm)\,$ in $\,\csg^\pm$ with
\begin{equation}
\phi^\pm\nh,\chi^\pm\in\csg^\pm\hs,\hskip7pt
|\phi^\pm|=|\chi^\pm|=1\hs,\hskip7pt\langle\phi^\pm\nh,\chi^\pm\rangle=0\hs,
\hskip7pt\phi^+\!\wedge\hs\chi^+=\hs\phi^-\!\wedge\hs\chi^-\nh.\label{bss}
\end{equation}
Every \spingeo\ $\,(\csg^+\!,\csg^-\!,\cka)\,$ gives rise to a {\em 
determinant mapping\/} $\,\,\det:\,\h(\csg^+\!,\csg^-)\to\bbC\,$ such that, 
for $\,A:\csg^+\to\csg^-\nh$, the operator 
$\,A^{\wedge2}:[\csg^+]^{\wedge2}\to[\csg^-]^{\wedge2}$ is the multiplication 
by $\,\hs\det A\,$ in the line 
$\,[\csg^+]^{\wedge2}=[\csg^-]^{\wedge2}=\cka$. Let us call 
$\,A\in\,\h(\csg^+,\csg^-)$ a {\em homothety\/} if 
$\,|A\phi\hs|\,=\,|A|\,|\phi\hs|\,$ for some $\,|A|\ge0$ and all 
$\,\phi\in\csg^+$, and set
\begin{equation}
\mathcal{V}\,=\,\{A\in\,\h(\csg^+\!,\csg^-):
\hskip3ptA\hskip7pt{\rm is\ a\ homothety\ and}\hskip5pt
\det A\in[\hs0,\infty)\}\,.\label{vee}
\end{equation}
\begin{rem}\label{revee}Any fixed $\,\phi^\pm\nh,\chi^\pm$ with (\ref{bss}) 
obviously make $\,\mathcal{V}\,$ correspond to the set of all $\,2\times2\,$ 
matrices appearing in (\ref{mlt}), with $\,a,b\in\bbC\hs$. Thus, 
$\,\mathcal{V}\,$ is a real vector space, $\,\dimr\mathcal{V}=4$, and 
$\,A\mapsto|A|\,$ is a Euclidean norm in $\,\mathcal{V}$, with 
$\,|A|^2=\hs\det A=|a|^2+|b|^2$ if $\,a,b\in\bbC\,$ represent $\,A\,$ as 
above. Also, $\,\mathcal{V}\,$ is canonically oriented, by (\ref{ori}), since 
every choice of bases $\,(\phi^\pm\nh,\chi^\pm)\,$ with (\ref{bss}) leads to 
an isomorphism $\,\mathcal{V}\to\bbC^2\nh$, and such 
$\,(\phi^\pm\nh,\chi^\pm)$ form an orbit of the connected Lie subgroup \hs 
Spin\cf\ \hs of $\,\hs\ug(2)\times\hs\ug(2)\,$ consisting of all 
$\,(\mathfrak A,\mathfrak B)\,$ with $\,\hs\det\mathfrak A=\det\mathfrak B$.
\end{rem}
\begin{lem}\label{spgeo}For any 
spin$^{\hskip.3pt{\rm c}}(4)$-ge\-om\-e\-try\/ $\,(\csg^+\!,\csg^-\!,\cka)\,$ 
there exist norm-pre\-serv\-ing real-linear isomorphic identifications\/ 
$\,\csg^\pm=\bbH\,$ and $\,\mathcal{V}=\bbH\hs$, which are also complex-linear 
for\/ $\,\csg^\pm$ and o\-ri\-en\-ta\-tion-pre\-serv\-ing for\/ 
$\,\mathcal{V}$, with\/ {\rm(\ref{neg})}, and which make the\/ {\em Clifford 
multiplication\/} $\,\mathcal{V}\times\csg^+\to\csg^-$, i.e., the evaluation 
pairing $\,(A,\phi)\mapsto A\phi$, appear as the quaternion multiplication. 
Specifically, such identifications are provided by declaring $\,\phi^\pm$ 
equal to\/ $\,1\in\bbH\,$ and\/ $\,\chi^\pm$ equal to\/ 
$\,\,{\bf j}\in\bbH\hs$, for any fixed\/ $\,(\phi^\pm\nh,\chi^\pm)\,$ with\/ 
{\rm(\ref{bss})}.
\end{lem}
This is clear from (\ref{mlt}) and Remark~\ref{revee}.

Insisting in Lemma~\ref{spgeo} that the identification $\,\mathcal{V}=\bbH\,$ 
be o\-ri\-en\-ta\-tion-pre\-serv\-ing might seem pedantic -- after all, the 
orientation of $\,\mathcal{V}\,$ can always be reversed. The point is, 
however, that given a \spst\ on a $\,4$\mfd\ $\,\y$, the canonical orientation 
of $\,\mathcal{V}\,$ leads to an orientation of $\,\y\,$ (see \S\ref{spct}). 
When the \spst\ comes from an \acst\ (along with a compatible metric, cf.\ 
\S\ref{acsu}), the latter orientation agrees with another canonical 
orientation of $\,\y$, provided by (\ref{ori}). Without this agreement, the 
generalization of formula~(\ref{tno}), obtained in \S\ref{trei}, might well 
read $\,\nu_f=T\x\,$ instead of $\,\nu_f=\hs\overline{T\x}$.

\section{Spin$^{\hskip.3pt{\rm c}}$-struc\-tures over \hs4-man\-i\-folds}
\label{spct}
By a {\em spin}$^{\hskip.3pt{\rm c}}${\it-struc\-ture\/} over a four\mfd\
$\,\y\,$ we mean a triple $\,(\sg^+\!,\sg^-\!,\ka)\,$ of complex vector
bundles of ranks $\,2$, $\,2\,$ and $\,1\,$ over $\,\y$, all endowed with 
Her\-mit\-i\-an fibre metrics, whose fibres form, at every point $\,y\in\y$, a 
\spingeo\ (Definition~\ref{defin}), varying with $\,y$, and such that the 
associated space $\,\mathcal{V}=\mathcal{V}_y$ with (\ref{vee}) has a fixed 
isomorphic identification with $\,T_y\y\,$ depending, along with the \spingeo\ 
itself, $\,C^\infty$-differentiably on $\,y$. (See \cite{lawson-michelsohn}, 
\cite{moore}.) In other words, we then have fixed norm-pre\-serv\-ing 
identifications 
$\,\ka\,=\,[\sg^+]^{\wedge2}=\,[\sg^-]^{\wedge2}$ and the {\em Clifford 
multiplication\/} $\,T\y\otimes\,\sg^+\to\sg^-$, which is a $\,C^\infty$ 
morphism  of real vector bundles acting, at every $\,y\in\y$, as a pairing 
$\,T_y\y\times\sg^+_y\ni(v,\phi)\mapsto v\phi\in\sg^-_y$ such that 
$\,\phi\mapsto v\phi\,$ equals, for any given $\,v$, a scalar $\,|v|\ge0\,$ 
times a norm-pre\-serv\-ing {\em unimodular\/} complex isomorphism 
$\,\sg^+_y\to\sg^-_y$. (Unimodularity makes sense here, since both $\,\sg^\pm$ 
have the same determinant bundle $\,\ka$.) This defines a Euclidean norm 
$\,v\mapsto|v|\,$ in $\,T_y\y$, and hence a Riemannian metric $\,g\,$ on 
$\,\y$. Also, $\,\y\,$ is canonically oriented (since so is the space 
$\,\mathcal{V}\,$ with (\ref{vee})). In other words, a \spst\ really lives 
over an {\em oriented Riemannian\/} $\,4$\mfd.

A special case is a {\em spin structure\/} $\,(\sg^+\!,\sg^-)\,$ for a 
$\,4$\mfd\ $\,\y$, that is, a \spst\ $\,(\sg^+\!,\sg^-\!,\ka)\,$ over $\,\y\,$ 
in which $\,\ka\,$ is the product line bundle $\,\y\times\,\bbC\,$ with the 
standard (constant) fibre metric. See \cite{milnor}, \cite{lawson-michelsohn}.

\section{Almost complex surfaces as 
spin$^{\hskip.3pt{\rm c}}$-man\-i\-folds}\label{acsu}
An {\em almost Her\-mit\-i\-an structure\/} $\,(J,g)\,$ on a manifold $\,\y\,$ 
consists of an \acst\ $\,J\,$ on $\,\y\,$ and a Riemannian metric $\,g\,$ on 
$\,\y$ compatible with $\,J$. Such $\,(J,g)\,$ may be viewed as a special case 
of a \spst\ (see \cite{lawson-michelsohn}, \cite{moore}), as described below 
in real dimension $\,4$.

More precisely, almost Hermitian structures on any $\,4$\mfd\ $\,\y\,$ are in 
a natural bijective correspondence with pairs formed by a \spst\ 
$\,(\sg^+\!,\sg^-\!,\ka)\,$ over $\,\y\,$ and a global unit $\,C^\infty$ 
section $\,\psi\,$ of $\,\sg^+$.

In fact, given such $\,(\sg^+\!,\sg^-\!,\ka)\,$ and $\,\psi$, the Clifford 
multiplication by $\,\psi\,$ identifies $\,T\y\,$ with $\,\sg^-$, thus 
introducing an \acst\ $\,J\,$ on $\,\y$, compatible with the Riemannian metric 
$\,g\,$ on $\,\y\,$ obtained as the real part of the Hermitian fibre metric 
$\,\langle\,,\rangle\,$ in $\,\sg^-$ (which, since $\,|\psi|=1$, is the same 
$\,g\,$ as in \S\ref{spct}). In addition, we have a natural isomorphic 
identification $\,\sg^+\nh=\iota\oplus\ka\,$ for the trivial complex line 
bundle $\,\iota=\y\times\bbC\hs$, where the subbundles 
$\,\hs\spanc\psi\,$ and $\,\psi^\perp$ of $\,\sg^+$ are identified with 
$\,\iota\,$ and $\,\ka\,$ via the norm-pre\-serv\-ing isomorphisms 
$\,\phi\mapsto\langle\psi,\phi\rangle\,$ and $\,\phi\mapsto\psi\wedge\phi$.

Conversely, any almost Her\-mit\-i\-an structure $\,(J,g)\,$ on $\,\y\,$ 
arises in this manner from $\,(\sg^+\!,\sg^-\!,\ka)$, where $\,\sg^-=T\y\,$ 
with the Hermitian fibre metric $\,\langle\,,\rangle\,$ whose real part is 
$\,g$, while $\,\ka=\detc\hs T\y\,$ and $\,\sg^+\nh=\iota\oplus\ka\,$ for 
$\,\iota=\y\times\bbC\,$ as above. This identifies $\,\ka\,$ with both 
$\,[\sg^\pm]^{\wedge2}$, as 
required. Also, $\,\sg^+$ is naturally a real vector 
subbundle of $\,\eta=\,\hr(T\y,T\y)$. Namely, the sections of $\,\iota\,$ 
provided by the constant functions $\,1,i:\y\to\bbC\,$ are to be identified 
with the sections $\,\,{\rm Id}\,\,$ and $\,J\,$ of $\,\eta$, while 
$\,\ka=\detc\hs T\y\,$ becomes a subbundle of $\,\eta\,$ if one lets 
$\,u\wedge v\in\ka_y$, for any $\,y\in\y$, operate on vectors $\,w\in T_y\y\,$ 
via $\,w\mapsto\langle v,w\rangle u-\langle u,w\rangle v$. The Clifford 
multiplication by $\,w\in T_y\y$, for $\,y\in \y$, then is the evaluation 
operator $\,\sg^+_y\ni A\mapsto Aw\in\sg^-_y$. (That $\,|Aw|=|A|\hs|w|\,$ is 
easily verified if one writes $\,A=(a,b\hs u\wedge v)\,$ with 
$\,a,b\in\bbC\,$ and $\,\langle\,,\rangle$-or\-tho\-nor\-mal vectors 
$\,u,v\in T_y\y$.) The distinguished unit section $\,\psi\,$ of $\,\sg^+$ is 
the constant function $\,1\,$ identified as above with a section of 
$\,\iota\subset\sg^+$, so that the Clifford multiplication by $\,\psi\,$ is a 
$\,\bbC$-linear bundle isomorphism $\,T\y\to T\y=\sg^-$. Therefore, these 
$\,(\sg^+\!,\sg^-\!,\ka)\,$ and $\,\psi\,$ in turn lead, via the construction 
of the preceding paragraph, to the original $\,(J,g)$.

\begin{rem}\label{reemb}For $\,n=2$, the Lie-group embedding 
$\,\hs\ug(n)\to$\hskip4ptSpin$^{\hskip.5pt{\rm c}}(2n)\,$ mentioned in 
\S\ref{intro} is given by 
$\,\mathfrak A\mapsto({\rm diag}\hs(1,\hs\det\mathfrak A),\hs\mathfrak A)$, in 
the notation of Remark~\ref{revee}. (See~\cite[p.~53]{moore}.) In fact, given 
an \acsu\ $\,\y\,$ and $\,y\in\y$, it is this homomorphism that renders 
equivariant the mapping which sends any orthonormal basis 
$\,(\phi^-\nh,\chi^-)$ of $\,\sg^-_y=T_y\y\,$ to the basis 
$\,(\phi^+\nh,\chi^+\nh,\phi^-\nh,\chi^-)$ of $\,\sg^+_y\oplus\hs\sg^-_y$ with 
(\ref{bss}), defined by $\,\phi^+=1\,$ and 
$\,\chi^+=\hs\phi^-\!\wedge\hs\chi^-\nh$, where $\,\sg^+_y=\bbC\oplus\ka_y$, 
and so $\,\phi^+\in\bbC\subset\sg^+_y$, while $\,\chi^+\in\ka_y\subset\sg^+_y$.

\end{rem}

\section{An explicit diffeomorphism 
$\,S^2\hskip-1.2pt\times S^2\,\approx\,G_2^+(\bbR^{\hskip-.3pt4})$}\label{xpdi}
A \spingeo\ (see Definition~\ref{defin}) gives rise to a specific 
diffeomorphic identification 
$\,S^2\hskip-1.2pt\times S^2\,\approx\,G_2^+(\bbR^{\hskip-.3pt4})$, with 
$\,G_2^+(\,)\,$ as in \S\ref{grop}.

Namely, given a \spingeo\ $\,(\csg^+\!,\csg^-\!,\cka)\,$ and any complex lines 
$\,\mathcal{L}\subset\csg^+$ and $\,\mathcal{L}'\subset\csg^-$ (cf.\ 
\S\ref{prel}), let $\,\mathcal{T}=\varPsi(\mathcal{L},\mathcal{L}')\,$ be the 
real subspace $\,\{A\in\mathcal{V}:A\mathcal{L}\subset\mathcal{L}'\}\,$ of the 
oriented real $\,4$-space $\,\mathcal{V}\,$ given by (\ref{vee}). Then
\begin{enumerate}
  \def\theenumi{{\rm(\alph{enumi}}}
\item $\mathcal{T}\,$ is naturally real-isomorphic to the complex line 
$\,\,\h(\mathcal{L},\mathcal{L}')$. Thus, $\,\dim\mathcal{T}=2\,$ and 
$\,\mathcal{T}\,$ is canonically oriented, i.e., 
$\,\mathcal{T}\in \gv$. 
\item $\overline{\mathcal{T}}{}^\perp
=\varPsi(\mathcal{L}^\perp\!,\mathcal{L}')$, where\/ 
$\,\overline{\mathcal{T}}{}^\perp\subset\mathcal{V}\,$ is the orthogonal 
complement of the oriented plane $\,\mathcal{T}$, endowed with the {\em 
nonstandard\/} orientation (cf.\ \S\ref{prel}).
\item Given $\,\psi\in\csg^+\smallsetminus\{0\}$, we have 
$\,\psi\in\mathcal{L}\,$ (or, $\,\psi\in\mathcal{L}^\perp$) if and only if the 
operator $\,\mathcal{T}\ni A\mapsto A\psi\in\csg^-$ is an o\-ri\-en\-ta\-tion 
pre\-serv\-ing (or, re\-vers\-ing) real-linear isomorphism onto a complex line 
in $\,\csg^-$.
\end{enumerate}
Finally, for $\,\varPsi(\mathcal{L},\mathcal{L}')=\mathcal{T}\,$ depending as 
above on complex lines $\,\mathcal{L},\mathcal{L}'\nh$,
\begin{equation}
\varPsi:P(\csg^+)\times P(\csg^-)\to \gv\hskip6pt
{\rm is\ a\ diffeomorphism,}\label{dif}
\end{equation}
with $\,\pw\,$ defined as in {\rm\S\ref{prel}} for $\,W=\csg^\pm$.

The remainder of this section is devoted to proving (a) -- (c) and 
(\ref{dif}). We set $\,\csg^\pm=\mathcal{V}=\bbH\hs$, as in Lemma~\ref{spgeo}, 
and choose $\,p,q\in\bbH\smallsetminus\{0\}\,$ with 
$\,\mathcal{L}=p\hs\bbC\hs$, $\,\hs\mathcal{L}'=q\hs\bbC\hs$. Then 
$\,\mathcal{T}=\{x\in\bbH:xp\in q\hs\bbC\}=q\hs\bbC\hs p^{-1}$, with 
$\,\bbC=\,\spanr(1,\hs{\bf i})\,\subset\,\bbH\hs$, so that 
(a) follows, and $\,(qp^{-1},\hs q\hs{\bf i}\hs p^{-1})$ is a 
pos\-i\-tive-o\-ri\-ented basis of $\,\mathcal{T}$.

Since the quaternion norm is multiplicative, the left or right multiplication 
by a nonzero quaternion is a homothety. Given 
$\,p,q\in\bbH\smallsetminus\{0\}$, both 
$\,(p,p\hs{\bf i},p\hs{\bf j},p\hs{\bf k})\,$ and 
$\,(q\hs p^{-1}\nh,q\hs{\bf i}\hs p^{-1}\nh,q\hs{\bf j}\hs p^{-1}\nh,
q\hs{\bf k}\hs p^{-1})\,$ thus are real-or\-thog\-o\-nal bases of $\,\bbH\hs$. 
For $\,\mathcal{L}=p\hs\bbC\,$ we now have 
$\,\mathcal{L}^\perp=p\hs{\bf j}\hs\bbC\,$ (as one sees using the first 
basis: $\,(p,\hs p\hs{\bf i}\hs)\,$ is an $\,\bbR$-basis of $\,\mathcal{L}$, 
and so $\,p\hs{\bf j}\hs\,$ and $\,-\hs p\hs{\bf k}=p\hs{\bf j}{\bf i}\hs\,$ 
must form an $\,\bbR$-basis of $\,\mathcal{L}^\perp$). Also, since 
$\,\mathcal{T}=\,\spanr(q\hs p^{-1}\nh,q\hs{\bf i}\hs p^{-1})$, orthogonality 
of the {\em second\/} basis implies that 
$\,\overline{\mathcal{T}}{}^\perp
=\,\spanr(q\hs{\bf j}\hs p^{-1}\nh,q\hs{\bf k}\hs p^{-1})$, i.e., 
$\,\overline{\mathcal{T}}{}^\perp=q\hs\bbC\hs(p\hs{\bf j})^{-1}
=\varPsi(\mathcal{L}^\perp,\mathcal{L}')$, which proves (b). (We write 
$\,\overline{\mathcal{T}}{}^\perp$ rather than $\,\mathcal{T}^\perp$, since 
the pos\-i\-tive-o\-ri\-ented basis $\,(qp^{-1}\nh,q\hs{\bf i}\hs p^{-1})\,$ 
of $\,\mathcal{T}\,$ and its analogue 
$\,(q(p{\bf j})^{-1}\nh,q\hs{\bf i}\hs(p{\bf j})^{-1})
=(-\hs q\hs{\bf j}\hs p^{-1}\nh,-\hs q\hs{\bf k}\hs p^{-1})$ for the 
orthogonal complement together form a {\em negative-oriented\/} basis of 
$\,\bbH=\mathcal{V}$, due to (\ref{neg}) and connectedness of 
$\,\bbH\smallsetminus\{0\}$.)

Furthermore, the $\,\bbC$-linear operator $\,\bbH\to\bbH\,$ of left 
multiplication by any fixed nonreal quaternion, being a matrix operator of the 
form (\ref{mlt}), must have a pair of conjugate nonreal eigenvalues. By 
Lemma~\ref{spgeo}, this applies to $\,\fo=A^{-1}\nh B:\csg^+\to\csg^+$ (or, 
$\,\fo=-AB^{-1}:\csg^-\to\csg^-$) 
whenever $\,A,B\in\mathcal{V}\,$ are linearly independent over $\,\bbR\hs$. 
For any $\,\mathcal{T}\in \gv$, let us now define 
$\,\varPsi'(\mathcal{T})\,$ to be the pair 
$\,(\mathcal{L},\mathcal{L}')\in P(\csg^+)\times P(\csg^-)\,$ obtained by 
fixing any pos\-i\-tive-o\-ri\-ented basis $\,(A,B)\,$ of $\,\mathcal{T}\,$ 
and choosing $\,\mathcal{L}\,$ (or, $\,\mathcal{L}'$) to be the eigenspace of 
$\,A^{-1}\nh B\,$ (or, $\,-AB^{-1}$) for the unique eigenvalue $\,z\in\bbC\,$ 
with $\,\,{\rm Im}\,z>0$. That $\,\varPsi'$ does not depend on the choice of 
the basis $\,(A,B)\,$ is easily seen, for $\,\mathcal{L}$, if we replace 
$\,(A,B)\,$ with $\,(B,-\hs A)\,$ or a basis of the form $\,(A,B')$, and, for 
$\,\mathcal{L}'\nh$, if instead of $\,(A,B)\,$ we use $\,(B,-A)\,$ or $\,(A',B)$. 
It is now immediate that $\,A\mathcal{L}\subset\mathcal{L}'$ and 
$\,B\mathcal{L}\subset\mathcal{L}'\nh$, i.e., 
$\,\varPsi(\mathcal{L},\mathcal{L}')=\mathcal{T}$, which shows that 
$\,\varPsi\circ\varPsi'=\,{\rm Id}$. Moreover, using Lemma~\ref{spgeo} as 
before, we get $\,\varPsi'\nh\circ\varPsi=\,{\rm Id}$. Namely, 
$\,\mathcal{T}=q\hs\bbC\hs p^{-1}$ has the pos\-i\-tive-o\-ri\-ented basis 
$\,(A,B)=(qap^{-1}\nh,qazp^{-1})$, with $\,a,z\in\bbC\hs$, $\,a\ne0\,$ and 
$\,\hs{\rm Im}\,z>0$. This makes $\,A^{-1}\nh B\,$ (or, $\,-AB^{-1}$) 
appear as the left quaternion multiplications by $\,pzp^{-1}$ (or, 
$\,-\hs qz^{-1}q^{-1}$), which has the eigenspace $\,\mathcal{L}=p\hs\bbC\,$
(or, $\,\mathcal{L}'=q\hs\bbC$) with the eigenvalue $\,z\,$ (or, $\,z^{-1}$). 
Now (\ref{dif}) follows, with $\,\varPsi'=\varPsi^{-1}$.

Finally, according to the above description of $\,\varPsi'=\varPsi^{-1}$ in 
terms of a pos\-i\-tive-o\-ri\-ented basis $\,(A,B)\,$ of $\,\mathcal{T}$, 
condition $\,\psi\in\mathcal{L}\,$ (or, $\,\psi\in\mathcal{L}^\perp$) is 
equivalent to $\,B\psi=zA\psi\,$ (or, $\,B\psi=\,\overline zA\psi$) for some 
$\,z\in\bbC\,$ with $\,\hs{\rm Im}\,z>0$. (The case of $\,\mathcal{L}^\perp$ 
follows from that of $\,\mathcal{L}\,$ since $\,\fo=A^{-1}\nh B\,$ is a 
homothety, cf.\ (\ref{vee}), and so $\,\mathcal{L}^\perp$ must be 
$\,\fo$-invariant as long as $\,\mathcal{L}\,$ is.) In other words, 
$\,\psi\in\mathcal{L}\,$ (or, $\,\psi\in\mathcal{L}^\perp$) if and only if 
$\,(A\psi,B\psi)\,$ is a positive (or, negative) oriented real basis of a 
complex line. This yields (c).

\section{Line bundles associated with a 
spin$^{\hskip.3pt{\rm c}}$-struc\-ture}\label{libu}
Given a \spst\ $\,(\sg^+\!,\sg^-\!,\ka)\,$ over a $\,4$\mfd\ $\,\y\,$ (see 
\S\ref{spct}), let $\,\la^+$ and $\,\la^-$ be the complex line bundles over 
the Grassmannian manifold $\,\gy\,$ (defined in \S\ref{grop}), whose fibres at 
any $\,(y,\mathcal{T})\,$ are the lines $\,\mathcal{L},\mathcal{L}'$ in 
$\,\sg^+_y,\hs\sg^-_y$ with 
$\,\varPsi(\mathcal{L},\mathcal{L}')=\mathcal{T}\,$ for $\,\varPsi\,$ as in 
(\ref{dif}), where $\,(\csg^+\!,\csg^-\!,\cka)=(\sg^+_y\!,\sg^-_y\!,\ka_y)\,$ 
and $\,\mathcal{V}$, given by (\ref{vee}), is identified with $\,T_y\y$. We 
then have
\begin{equation}
{\rm a)}\hskip13pt\pi^*\sg^\pm\,=\,\,\la^\pm\oplus\,\mu^\pm\nh,\hskip43pt
{\rm b)}\hskip13pt\pi^*\ka\,=\,\la^\pm\otimes\,\mu^\pm\nh,\label{psg}
\end{equation}
for the bundle projection $\,\pi:\gy\to\y\,$ and the line bundles $\,\mu^\pm$ 
which are the orthogonal complements of $\,\la^\pm$ in $\,\pi^*\sg^\pm$. (The 
natural isomorphic identification in (\ref{psg}\hs-b) is obvious from 
(\ref{psg}\hs-a) as $\,\ka=[\sg^\pm]^{\wedge2}$, cf.\ \S\ref{spct}.)

We wish to emphasize that even when the \spst\ $\,(\sg^+\!,\sg^-\!,\ka)\,$ in 
question arises from an almost Hermitian structure $\,(J,g)$, as in 
\S\ref{acsu}, $\,J\,$ and $\,g\,$ do not seem to provide any shortcuts for 
defining $\,\la^\pm$ and $\,\mu^\pm$. If anything, they are rather an 
impediment, unless one simply ignores them; this amounts to ignoring both the 
decomposition $\,\sg^+\nh=\iota\oplus\ka\,$ and the presence of a 
distinguished unit section $\,\psi\,$ of $\,\sg^+$ and, in effect, treating 
$\,(\sg^+\!,\sg^-\!,\ka)\,$ as if it were just any \spst, with no additional 
features.

For any oriented $\,4$\mfd\ $\,\y\,$ and $\,\gy,\hs\pi\,$ as above, the 
tangent bundle $\,\tgy\,$ admits specific oriented real-plane subbundles 
$\,\tau,\nu\,$ such that $\,\pi^*T\y=\tau\oplus\nu\,\subset\,\tgy$. Namely, 
$\,\tau\,$ is the {\em tautological bundle\/} with the fibre 
$\,\mathcal{T}\subset T_y\y\,$ over any $\,(y,\mathcal{T})$. The other summand 
$\,\nu$, which might be called the {\em tautonormal bundle}, is obtained by 
choosing a Riemannian metric on $\,\y\,$ and setting 
$\,\nu=\tau^\perp\subset\pi^*T\y$.

If $\,\y\,$ now happens to be the canonically-oriented base manifold of a 
\spst\ $\,(\sg^+\!,\sg^-\!,\ka)$, the definitions of $\,\la^\pm\nh,\mu^\pm$ 
and (a), (b) in \S\ref{xpdi} give the following natural isomorphic 
identifications of oriented real-plane bundles:
\begin{equation}
\tau\,=\,\hs\hc(\la^+,\la^-)\,,\hskip30pt
\overline\nu\,=\,\hs\hc(\mu^+,\la^-)\,,\label{tnu}
\end{equation}
with $\,\overline\nu\,$ standing for $\,\nu\,$ with the reversed orientation. 
We may treat $\,\tau,\nu\,$ as complex line bundles, using their orientations 
and the fibre metrics on them induced by a fixed Riemannian metric on $\,\y$, 
along with Remark~\ref{relin} (or, just (\ref{tnu})). Then 
$\,\overline\nu=\nu^*$ and $\,\hs\hc(\overline\nu,\tau)=\nu\otimes\tau$, so 
that, by (\ref{tnu}), 
$\,\nu\otimes\tau=\mu^+\otimes\overline{\la^-}\otimes\la^-\otimes
\overline{\la^+}=\mu^+\otimes\overline{\la^+}$, i.e., a natural isomorphism 
also exists between $\,\hs\hc(\overline\nu,\tau)\,$ and 
$\,\hs\hc(\la^+,\mu^+)$.

Since $\,T\x=F^*\tau\,$ and $\,\nu_f=F^*\nu$, where $\,\nu_f$ and 
$\,\hf:\x\to\gy\,$ are the normal bundle and the Gauss mapping of any given 
immersion $\,f:\x\to\y\,$ of an oriented real surface $\,\x\,$ (see 
\S\ref{grop}), this leads to an {\em adjunction formula}, which equates 
$\,\hs\hc(\overline{\nu_f},T\x)\,$ with the $\,F$-pull\-back of 
$\,\hs\hc(\la^+\nh,\mu^+)$.

\section{The mapping \hs{\rm(\ref{dif})}\hs\ and antilinearity}
\label{tvcs}
Let $\,\mathcal{V}\,$ be the real $\,4$-space with (\ref{vee}) for a \spingeo\ 
$\,(\csg^+\!,\csg^-\!,\cka)$, and let $\,\varPsi\,$ be the diffeomorphism 
appearing in (\ref{dif}). Then, for any $\,\mathcal{L}'\in P(\csg^-)$,
\begin{equation}
{\rm the\ differential\ of}\hskip9pt
P(\csg^+)\ni\mathcal{L}\mapsto\varPsi(\mathcal{L},\mathcal{L}')\in 
\gv\hskip9pt{\rm is\ antilinear}\label{inc}
\end{equation}
at every point of $\,P(\csg^+)$, antilinearity referring to the standard 
complex structure of the complex projective line $\,P(\csg^+)\,$ and the 
\acst\ on $\,\gv\,$ associated, as in \S\ref{grop}, with the Euclidean inner 
product in $\,\mathcal{V}\,$ mentioned in  Remark~\ref{revee}.

In fact, let $\,\csg^\pm=\mathcal{V}=\bbH\,$ as in Lemma~\ref{spgeo}. 
As stated in to the two paragraphs following (\ref{dif}), for any 
given $\,\mathcal{L}\in P(\csg^+)\,$ and $\,\mathcal{L}'\in P(\csg^-)\,$ we 
then have $\,\mathcal{L}=p\hs\bbC\hs$, 
$\,\mathcal{L}^\perp=p\hs{\bf j}\hs\bbC\,$ and $\,\hs\mathcal{L}'=q\hs\bbC\,$ 
with some unit quaternions $\,p,q$, and so 
$\,(q\hs\overline p,\hs q\hs{\bf i}\hs\overline p\hs)\,$ is a 
pos\-i\-tive-o\-ri\-ented basis of the plane 
$\,\mathcal{T}=\varPsi(\mathcal{L},\mathcal{L}')$, i.e., 
$\,\mathcal{T}=q\hs\bbC\hs\overline p$, while 
$\,\mathcal{T}^\perp=q\hs\bbC\hs{\bf j}\hs\overline p$. (Note that 
$\,p^{-1}=\hs\overline p\,$ since $\,|p|=1$.) As usual, 
$\,\bbC=\,\spanr(1,\hs{\bf i})\,\subset\,\bbH\hs$. Thus, we may write 
$\,\mathcal{T}=\bbC\hs$, identifying $\,z\in\bbC\,$ with $\,qz\overline p$. 
The above basis of $\,\mathcal{T}\,$ then becomes $\,(1,i)$, so that 
$\,\bbC\,$ with its usual structure is precisely the complex line formed by 
the oriented Euclidean plane $\,\mathcal{T}\,$ (see Remark~\ref{relin}).

By (\ref{tlp}), 
$\,T_{\mathcal{L}}[P(\csg^+)]=\hs\hc(\mathcal{L},\mathcal{L}^\perp)
=\hs\hc(p\hs\bbC,p\hs{\bf j}\hs\bbC)=\bbC\hs$, where each $\,c\in\bbC\,$ is 
identified with the operator $\,p\hs\bbC\to p\hs{\bf j}\hs\bbC\,$ sending 
$\,p\,$ to $\,p\hs{\bf j}\hs c$. Similarly, the description of 
$\,T_L[G_m(W)]\,$ in \S\ref{prel} identifies the tangent space of 
$\,\gv\,$ at $\,\mathcal{T}\,$ with 
$\,\hs\hr(\mathcal{T},\mathcal{V}/\mathcal{T})
=\hs\hr(\mathcal{T},\mathcal{T}^\perp)$. We will now show that the 
differential of the mapping (\ref{inc}) at $\,\mathcal{L}\,$ sends every 
$\,c\in\bbC=T_{\mathcal{L}}[P(\csg^+)]\,$ to the operator 
$\,A_c:\bbC\to\mathcal{T}^\perp$ given by 
$\,A_cz=-\hs qz\hs\overline c\hs{\bf j}\hs\overline p\,$ (with 
$\,\mathcal{T}=\bbC\,$ as before). Since $\,A_{ic}z=-\hs A_c(iz)$, (\ref{inc}) 
will follow.

To this end, let us replace $\,p\,$ by a $\,C^1$ function of a real parameter 
$\,t$, valued in unit quaternions, equal at $\,t=0\,$ to the original $\,p$, 
and such that $\,dp/dt\,$ at $\,t=0\,$ equals $\,p{\bf j}c$. Then 
$\,d\hs{\bf w}/dt\,$ at $\,t=0\,$ for the $\,t$-dependent basis 
$\,{\bf w}\hs=(q\hs\overline p,\hs q\hs{\bf i}\hs\overline p)\,$ of the 
$\,t$-dependent plane $\,\mathcal{T}=q\hs\bbC\hs\overline p\,$ equals 
$\,{\bf v}\hs=(A_c1,A_ci)$, with $\,A_c$ defined above, which (cf.\ 
Remark~\ref{regra}) proves our claim about $\,A_c$.
\begin{rem}\label{rebih}A virtually identical argument shows that the 
differential of 
$\,P(\csg^-)\ni\mathcal{L}'\mapsto\varPsi(\mathcal{L},\mathcal{L}')\in\gv\,$ 
is {\em complex-linear\/} at every point of $\,P(\csg^-)$, for any fixed 
$\,\mathcal{L}'\in P(\csg^-)$. (Note that in 
$\,\mathcal{T}=q\hs\bbC\hs\overline p\,$ conjugation is applied to $\,p$, but 
not to $\,q$.) Along with (\ref{inc}) this implies that the \acst\ of 
$\,\gv\,$ is integrable and (\ref{dif}) is a biholomorphism, provided that the 
factor $\,P(\csg^+)\,$ in $\,P(\csg^+)\times P(\csg^-)\,$ carries the {\em 
conjugate\/} of its standard complex structure.

Note the well-known fact that, for a Euclidean space $\,W\,$ of any dimension, 
the \acst\ on $\,G_2^+(W)\,$ described in \S\ref{grop} is integrable.
\end{rem}

\section{Transversality}\label{trsv}
Given a $\,C^1$ section $\,\phi\,$ of a vector bundle $\,\eta\,$ over a 
manifold $\,N\hs$ and a point $\,\xi\in N\,$ with $\,\phi(\xi)=0$,
\begin{enumerate}
  \def\theenumi{{\rm(\alph{enumi}}}
\item $\phi\,$ is transverse at $\,\xi\,$ to the zero section of 
$\,\eta\,$ if and only if, for the fibre-valued function 
$\,\phi\hs':U\to\bbF\,$ representing $\,\phi\,$ in some, or any, local 
trivialization of $\,\eta\,$ over a neighborhood $\,\,U\,$ of $\,\xi$, 
the differential $\,d\phi_\xi':T_\xi N\to\bbF\,$ is surjective.
\item $\phi\,$ is transverse at $\,\xi\,$ to the zero section of 
$\,\eta\,$ whenever this is the case for $\,\eta,\phi\,$ restricted to a 
submanifold of $\hs N$ containing $\,\xi$.
\end{enumerate}
This is a trivial exercise; (b) easily follows from (a).

Let $\,\pw\,$ now be a real, complex or quaternion projective line obtained as 
in \S\ref{prel} from a plane $\,W\,$ over a field $\,\bbK\,$ (one of 
$\,\bbR\hs,\bbC\hs,\bbH$), and let $\,\eta\,$ be the tautological 
$\,\bbK$-line bundle over $\,\pw$, with the fibre $\,L\,$ at any $\,L\in \pw$. 
Choosing, in addition, a $\,\bbK$-val\-ued sesqui\-linear inner product 
$\,\langle\,,\rangle\,$ in $\,W\,$ and a nonzero vector $\,u\in W$, we now 
also define a $\,\bbK$-line bundle $\,\zeta\,$ over $\,\pw\,$ along with 
$\,C^\infty$ sections $\,\phi\,$ in $\,\eta\,$ and $\,\chi\,$ in $\,\zeta$. 
Namely, $\,\zeta=\eta^\perp$ is the orthogonal complement of $\,\eta\,$ 
treated as a subbundle of the product bundle $\,\theta=\pw\times W$, while 
$\,\phi,\chi\,$ are the $\,\eta\,$ and $\,\zeta\,$ components of $\,u$, which 
is a constant section of $\,\theta$, relative to the decomposition 
$\,\theta=\eta\oplus\zeta$. Then
\begin{enumerate}
  \def\theenumi{{\rm(\roman{enumi}}}
\item Either of $\,\phi,\chi\,$ has just one zero, at which it is transverse 
to the respective zero section in $\,\eta\,$ or $\,\zeta$.
\item If $\,\bbK=\bbC\hs$, then, for some $\,\phi\hs',\chi\hs':U\to\bbC\,$ 
representing $\,\phi,\chi\,$ as in (a) above, the differential of 
$\,\phi\hs'$ (or, $\,\chi\hs'$) at the zero in question is an antilinear (or, 
complex-linear) isomorphism.
\end{enumerate}
In fact, we may fix $\,v\in W\,$ with $\,\langle u,v\rangle=0\,$ and 
$\,|v|=|u|$. The unique zero of $\,\phi\,$ (or, $\,\chi$) is 
$\,u^\perp=\bbK v\in \pw\,$ (or, $\,v^\perp=\bbK u\in\pw$); let us agree to 
identify its neighborhood $\,\y'=\pw\smallsetminus\{v^\perp\}\,$ (or, 
$\,\y'=\pw\smallsetminus\{u^\perp\}$) with $\,\bbK\hs$, using the 
diffeomorphism $\,\bbK\to\y'$ given by $\,z\mapsto\bbK(zu+v)\,$ (or, 
$\,z\mapsto\bbK(zv+u)$), cf.\ Remark~\ref{redvt}. Under this identification, 
a local trivializing section of $\,\eta\,$ (or, $\,\zeta$) is defined on 
$\,\y'=\bbK\,$ by $\,z\mapsto zu+v\,$ (or, $\,z\mapsto\overline zu+v$) and, 
evaluating $\,\langle\,,\rangle$-or\-thog\-o\-nal projections of $\,u\,$ onto 
the directions of $\,zu+v\,$ and $\,\overline zu+v\,$ in $\,W$, we see that 
the $\,\bbK$-val\-ued functions $\,\phi\hs',\chi\hs'$ corresponding to 
$\,\phi,\chi\,$ as in (a) above are $\,z\mapsto\overline z\,$ and 
$\,z\mapsto z\,$ divided, in both cases, by $\,(|z|^2+1)|u|^2$. Their 
differentials at the respective zeros (i.e., at $\,z=0$) are 
$\,z\mapsto\overline z/|u|^2\,$ and $\,z\mapsto z/|u|^2$. This gives (ii) and, 
combined with (a) above, proves the transversality claim in (i).
\begin{rem}\label{reori}Suppose that $\,K,K'$ are $\,2$\diml\ submanifolds of 
a real $\,4$\mfd\ $\,\y$, both containing some given point $\,\xi\in N$, and 
$\,K\,$ is the set of zeros of a $\,C^\infty$ section $\,\phi\,$ of a complex 
line bundle $\,\la\,$ over $\,N$, transverse to the zero section. Let 
$\,\ap,\be,\ga\,$ be the differentials at $\,\xi\,$ of the inclusion mappings 
$\,K\to N$, $\,K'\to N\,$ and, respectively, of the restriction to 
$\,K'\cap\,U$ of a function  $\,\phi\hs':U\to\bbC\,$ that represents 
$\,\phi\,$ in some local trivialization of $\,\la\,$ having a domain $\,\,U\,$ 
with $\,\xi\in U$. Finally, let $\,N,K\,$ and $\,K'$ all carry some fixed 
\acst s, so that they are all oriented via (\ref{ori}).

If $\,\ga\,$ is injective and $\,\be\,$ is antilinear, while $\,\ap,\ga\,$ 
are both complex-linear or both antilinear, then the orientation of $\,K\,$ 
defined by (\ref{ori}) is the opposite of the orientation that 
$\,K\subset N\,$ acquires by being the zero set of $\,\phi$, cf.\ \S\ref{prel}.

In fact, as $\,\phi\hs'=0\,$ on $\,K\cap\,U$, injectivity of $\,\ga\,$ gives 
$\,T_\xi N=T_\xi K\oplus T_\xi K'$. The complex-plane orientation of 
$\,T_\xi N\,$ (given by (\ref{ori})) is the direct sum of the complex-line 
orientation in $\,T_\xi K\,$ and the opposite of the orientation in 
$\,T_\xi K'$ pulled back 
by $\,\ga\,$ from $\,\bbC\hs$. In fact, our assumptions state that the complex 
structure prescribed in $\,T_\xi K$, and the one in $\,T_\xi K'$ making 
$\,\ga\,$ antilinear, are either both identical with, or both conjugates of, 
the complex structures that $\,T_\xi K,\hs T_\xi K'$ inherit from 
$\,T_\xi N\,$ by being its complex subspaces.
\end{rem}
\section{Poincar\'e duals 
of\hskip6pt$\pep$\hskip6ptand\hskip6pt$\pem$\hskip6ptin\hskip6pt$\gy$}
\label{podu}
Let $\,K\,$ be a compact, oriented, not necessarily connected, codimension-two 
submanifold of an oriented even\diml\ real manifold $\,N$, and let $\,\la\,$ 
be a complex line bundle over $\,N$. We will say that $\,\la\,$ is {\em 
Poincar\'e-dual\/} to $\,K$ if $\,K\,$ is the canonically-oriented set of 
zeros (cf.\ \S\ref{prel}) of some $\,C^\infty$ section of $\,\la\,$ transverse 
to the zero section.

This implies that $\,c_1(\la)\in H^2(N,\hs\bbZ)\,$ corresponds to 
$\,[K]\in H_{n-2}(N,\hs\bbZ)\,$ under the Poincar\'e duality, as 
$\,\int_\Delta c_1(\la)=[K]\cdot[\Delta]\,$ for any $\,2$-cycle $\,\Delta\,$ 
in $\,N$, due to the Poincar\'e index formula for $\,\la\,$ restricted to 
$\,\Delta$.
\begin{thm}\label{thuno}Let\/ $\,\y\,$ be an almost Hermitian surface treated 
as a\/ $\,4$\mfd\ along with a \spst\/ $\,(\sg^+\!,\sg^-\!,\ka)\,$ and a fixed 
global unit\/ $\,C^\infty$ section\/ $\,\psi\,$ of\/ $\,\sg^+$, cf.\ 
{\rm\S\ref{acsu}}. The line bundles\/ $\,\la^+$ and\/ $\,\mu^+$ over\/ 
$\,\gy$, introduced in\/ {\rm\S\ref{libu}}, then are Poincar\'e-dual to the 
oriented\/ $\,6$\diml\ \sbms\/ $\,\pem\,$ and\/ $\,\pep\,$ of\/ $\,\gy$, 
defined at the end of\/ {\rm\S\ref{grop}}. 

More precisely, the\/ $\,\la^+$ and $\,\mu^+$ components\/ $\,\phi,\chi\,$ 
of\/ $\,\pi^*\psi\,$ relative to the decomposition\/ {\rm(\ref{psg}\hs-a)} are 
transverse to the zero sections in\/ $\,\la^+$ and\/ $\,\mu^+\nh$, and their 
respective oriented manifolds of zeros are\/ $\,\pem\,$ and\/ $\,\pep$.
\end{thm}
\begin{proof}Let us fix $\,y\in\y\,$ and consider the \spingeo\ 
$\,(\csg^+\!,\csg^-\!,\cka)=(\sg^+_y,\sg^-_y,\ka)$, with the identification 
$\,\mathcal{V}=T_y\y\,$ (see \S\ref{spct}). The fibre $\,G_2^+(T_y\y)\,$ of 
$\,\gy\,$ over $\,y\,$ thus is identified with $\,\gv$. That 
the respective {\it sets\/} of zeros of $\,\phi,\chi\,$ in $\,\gy\,$ are 
$\,\pem\,$ and $\,\pep\,$ is now clear from (c) in \S\ref{xpdi} along with the 
fact that the Clifford multiplication by $\,\psi\,$ is an isomorphism 
$\,T\y\to\sg^-$ of {\em complex\/} vector bundles (see \S\ref{acsu}), and 
the definitions of $\,\la^+,\mu^+$ in \S\ref{libu}.

The same definitions clearly imply that the pull\-backs of $\,\la^+,\mu^+$ 
under (\ref{inc}), for any fixed $\,\mathcal{L}'$, are the bundles 
$\,\eta,\zeta\,$ over $\,\pw$, for $\,W=\csg^+$, described in \S\ref{trsv}, 
while $\,\phi\,$ and $\,\chi$, pulled back to $\,P(\csg^+)\,$ via (\ref{inc}), 
coincide with the sections $\,\phi\,$ and $\,\chi\,$ defined in \S\ref{trsv} 
for $\,u=\psi(y)$. Our transversality claim thus is obvious from (i) in 
\S\ref{trsv} combined with (b) in \S\ref{trsv} for the submanifold of 
$\,N=\gy\,$ obtained as the image of (\ref{inc}).

What remains to be shown is that the orientations of $\,\pem\,$ and $\,\pep\,$ 
defined at the end of \S\ref{grop} coincide with the orientations which they 
acquire by being the zero sets of transverse sections in complex line bundles 
over the oriented $\,8$\diml\ manifold $\,\gy\,$ (cf.\ \S\ref{prel}). 

However, $\,\pem,\pep\,$ and $\,\gy\,$ are bundle spaces over $\,\y$, and 
their orientations described in \S\ref{grop} are opposite to the direct sums 
of the 
orientations of the base $\,\y\,$ and those of the fibres $\,\pmv,\hs\pv\,$ 
or $\,\gv$. All four orientations are induced by 
\acst s via (\ref{ori}). It therefore suffices to establish agreement between 
the standard orientations of the complex projective lines $\,\pmv,\hs\pv$, 
and the orientations of $\,\pmv,\hs\pv\,$ as submanifolds of the canonically 
oriented \acsu\ $\,\gv\,$ which are the zero sets of 
$\,\phi,\chi\,$ restricted to the fibre $\,\gv$.

Such agreement is in turn immediate from Remark~\ref{reori} applied to 
$\,N=\gv$, $\,\xi=\mathcal{T}\,$ for any fixed $\,\mathcal{T}\in\pmv\,$ (or, 
$\,\mathcal{T}\in\pv$), $\,K'$ which is the image of (\ref{inc}) with 
$\,\mathcal{L}'$ chosen so that 
$\,\mathcal{T}=\varPsi(\mathcal{L},\mathcal{L}')\,$ for some $\,\mathcal{L}$, 
along with $\,K=\pmv\,$ (or, $\,K=\pv$) and with $\,\phi\,$ which is the 
restriction to the fibre $\,\gv\,$ of our $\,\phi\,$ (or, of our 
$\,\chi$). That the assumptions listed in Remark~\ref{reori} are all satisfied 
is clear from (ii) in \S\ref{trsv}, (\ref{dif}), (\ref{inc}) and 
Remark~\ref{reppg}. This completes the proof.
\end{proof}
\begin{cor}\label{corol}For\/ any almost Hermitian structure\/ $\,(J,g)\,$ on 
a\/ $\,4$\mfd\/ $\,\y$, the line bundle\/ $\,\pi^*\ka\,$ over\/ $\,\gy$, 
with\/ $\,\ka=\detc\hs T\y\,$ defined as in\/ {\rm\S\ref{cpts}}, is 
Poincar\'e-dual to the union\/ $\,\pep\cup\pem\subset\gy$. Thus, 
{\rm(\ref{gss})} holds for every immersion\/ $\,f:\x\to\y\,$ of a closed 
oriented real surface\/ $\,\x$.
\end{cor}
In fact, for $\,\phi,\chi\,$ as in Theorem~\ref{thuno}, $\,\phi\otimes\chi\,$ 
is transverse to the zero section in $\,\pi^*\ka=\la^+\otimes\,\mu^+$ (cf.\ 
(\ref{psg}\hs-b)), by (a) in \S\ref{trsv} with 
$\,\pep\cap\pem=$\hskip4.5pt\emp\hs.

One could also prove Corollary~\ref{corol} without mentioning 
$\,(\sg^+\!,\sg^-\!,\ka)\,$ and $\,\psi\,$ at all. Instead, one might use the 
section of $\,\pi^*\ka$, for $\,\ka=\detc\hs T\y$, which assigns to 
$\,(y,\mathcal{T})\in\gy\,$ the complex 
exterior product $\,v\wedge w\,$ for any pos\-i\-tive-o\-ri\-ented 
$\,g$-or\-tho\-nor\-mal basis $\,(v,w)\,$ of the plane 
$\,\mathcal{T}\subset T_y\y$. This section is transverse to the zero section 
in $\,\pi^*\ka$, as it equals $\,-\hs2i\hs\phi\otimes\chi\,$ (by 
Remark~\ref{rewed} below); however, having to establish its transversality 
directly would make such a proof quite tedious.

\section{Spin$^{\hskip.3pt{\rm c}}$-man\-i\-folds with Higgs fields
}\label{hggs}
We borrow the term 'Higgs field' from theoretical particle physics, where it 
stands for a distinguished section of a specific vector bundle over the 
spacetime manifold. Such a section leads to spontaneous symmetry breaking, 
that is, the reduction of the original structure group to a lower\diml\ Lie 
subgroup. (See, for instance, \S11.5 of \cite{derdzinski}.)

This is analogous to the situation in \S\ref{acsu}, where by fixing a 
global unit $\,C^\infty$ section $\,\psi\,$ of $\,\sg^+$ we reduce a \spst\ 
$\,(\sg^+\!,\sg^-\!,\ka)$, having the $\,7$\diml\ structure group \hs 
Spin\cf\ (see Remark~\ref{revee}), to an almost Her\-mit\-i\-an structure, 
with the $\,4$\diml\ group $\,\hs\ug(2)$. Allowing $\,\psi\,$ to have zeros, 
we arrive at the following generalization of almost Hermitian structures. 

Let $\,(\sg^+\!,\sg^-\!,\ka)\,$ be a \spst\ over a $\,4$\mfd\ $\,\y$. By a 
{\em Higgs field\/} for $\,(\sg^+\!,\sg^-\!,\ka)\,$ we mean a $\,C^\infty$ 
section $\,\psi\,$ of $\,\sg^+$ defined only up to multiplication by positive 
$\,C^\infty$ functions on $\,\y\,$ and transverse to the zero section. The 
product of such a section $\,\psi\,$ and a $\,C^\infty$ function 
$\,\y\to(0,\infty)\,$ then represents the same Higgs field as $\,\psi\,$ does; 
for simplicity, however, $\,\psi\,$ itself will also be called a 'Higgs field'.

For instance, any almost Hermitian structure on a $\,4$\mfd\ $\,\y\,$ gives 
rise to the \spst\ $\,(\sg^+\!,\sg^-\!,\ka)$, defined in \S\ref{acsu}, along 
with the Higgs field without zeros represented by the unit $\,C^\infty$ 
section $\,\psi\,$ of $\,\sg^+$ described near the end of \S\ref{acsu}. On the 
other hand, every \spst\ $\,(\sg^+\!,\sg^-\!,\ka)\,$ over a compact $\,4$\mfd\ 
admits a Higgs field $\,\psi$, which may be chosen arbitrarily $\,C^1$-close 
to any given $\,C^\infty$ section of $\,\sg^+$.
\begin{rem}\label{reres}Let $\,\zr\subset\y\,$ be the discrete set of 
all zeros of a Higgs field $\,\psi\,$ for a \spst\ $\,(\sg^+\!,\sg^-\!,\ka)\,$ 
over a $\,4$\mfd\ $\,\y$. The open submanifold  
$\,\y'=\y\nh\smallsetminus\zr\,$ then carries the {\em residual\/} \acst\ 
$\,J$, determined as in \S\ref{acsu} by $\,(\sg^+\!,\sg^-\!,\ka)\,$ restricted 
to $\,\y'$ and $\,\psi/|\psi|$. 
\end{rem}
\begin{lem}\label{trhgs}Given a Higgs field\/ $\,\psi\,$ for a \spst\ 
$\,(\sg^+\!,\sg^-\!,\ka)\,$ over a\/ $\,4$\mfd\ $\,\y$, let\/ 
$\,\phi,\chi\,$ be the\/ $\,\la^+$ and\/ $\,\mu^+$ components of 
the section\/ $\,\pi^*\psi\,$ of\/ $\,\pi^*\sg^+$ relative to the 
decomposition\/ {\rm(\ref{psg}\hs-a)}. Then\/ $\,\phi$ and\/ $\,\chi\,$ 
are transverse to the zero sections in\/ $\,\la^+$ and, respectively, 
$\,\mu^+$. 
\end{lem}
\begin{proof}Transversality at points of the open dense set 
$\,\pi^{-1}(\y')\subset\gy$, with $\,\y'=\y\nh\smallsetminus\zr\,$ as in 
Remark~\ref{reres}, is obvious from the final clause of Theorem~\ref{thuno}, 
applied to $\,(\sg^+\!,\sg^-\!,\ka)\,$ (and $\,\psi$) restricted to 
$\,\y'\nh$. Note that we may assume that $\,|\psi|=1\,$ on $\,\y'$ by suitably 
rescaling both fibre metrics of $\,\sg^\pm$, namely, dividing them by 
$\,|\psi|^2$. 

Given fibre-valued functions $\,\phi\hs'\nh,\chi\hs'$ that represent 
$\,\phi,\chi\,$ in local trivializations of $\,\la^+$ and $\,\mu^+$, it is 
clear that $\,(\phi\hs'\nh,\chi\hs')\,$ will similarly represent 
$\,\pi^*\psi=\phi+\chi\,$ in the corresponding direct-sum local 
trivialization of $\,\pi^*\sg^+=\la^+\nh\oplus\mu^+\nh$. Our transversality 
assertion at points of $\,\pi^{-1}(\zr)\,$ is now immediate from 
(a) in \S\ref{trsv}, which completes the proof.
\end{proof}
At the end of \S\ref{grop}\ we associated with any \acsu\ $\,\y\,$ a pair 
$\,\pep,\pem\,$ of oriented $\,6$\diml\ \sbms\ of the Grassmannian manifold 
$\,\gy$. Lemma~\ref{trhgs} provides a natural generalization of that 
construction to the case where $\,\y\,$ is a $\,4$\mfd\ carrying a \spst\ 
$\,(\sg^+\!,\sg^-\!,\ka)\,$ with a fixed Higgs 
field $\,\psi$. Namely, we declare $\,\pem\,$ and $\,\pep\,$ to be the zero 
sets in $\,\gy\,$ of the sections $\,\phi,\chi\,$ of $\,\la^+$ and $\,\mu^+$ 
obtained as the components of $\,\pi^*\psi\,$ relative to the decomposition 
(\ref{psg}\hs-a). In view of Lemma~\ref{trhgs}, $\,\pep,\pem\,$ are, again, 
not-necessarily-connected, oriented $\,6$\diml\ \sbms\ of $\,\gy$. The 
conclusion of Theorem~\ref{thuno} also remains valid, although now it is 
nothing else than the definition of  $\,\pem\,$ and $\,\pep$. 

Given $\,\y,(\sg^+\!,\sg^-\!,\ka),\psi\,$ as above and an immersion 
$\,f:\x\to\y\,$ of an oriented real surface $\,\x$, we call $\,x\in\x\,$ a 
{\em complex point\/} of $\,f\,$ if it is a $\,\pep${\it\hskip-1.4pt-com\-plex 
point\/} or a $\,\pem${\it\hskip-1.4pt-com\-plex point}, in the sense that 
$\,\hf(x)\in\pep\,$ or, respectively, $\,\hf(x)\in\pem$, where 
$\,\hf:\x\to\gy\,$ is the Gauss mapping of $\,f\,$ (see end of \S\ref{grop}). 
This generalizes of the definitions given in \S\ref{cpts}, and will in turn 
allow us to define, in \S\ref{trei} and \S\ref{pshi}, what it means for such 
an immersion $\,f\,$ to be \tr\ or \psh\ (which, again, is a straightforward 
extension from the almost-complex case).
\begin{rem}\label{recpp}For $\,\y,(\sg^+\!,\sg^-\!,\ka),\psi\,$ and $\,\zr\,$ 
as in Remark~\ref{reres}, and for any immersion $\,f:\x\to\y\,$ of an oriented 
real surface $\,\x$, every point $\,x\in\x\,$ with $\,f(x)\in\zr\,$ is a 
complex point of $\,f$. More precisely, every oriented plane in 
$\,T_{f(x)}\y\,$ then lies in both $\,\pep$ and $\,\pem$.
\end{rem}

\section{Totally real immersions of oriented surfaces}\label{trei}
Given a $\,4$\mfd\ $\,\y\,$ and a \spst\ $\,(\sg^+\!,\sg^-\!,\ka)\,$ over 
$\,\y\,$ with a fixed Higgs field $\,\psi\,$ (see \S\ref{hggs}), we will say 
that an immersion $\,f:\x\to\y\,$ of an oriented real surface $\,\x\,$ is {\em 
totally real\/} if $\,f\,$ has no complex points (defined as in \S\ref{hggs}), 
that is, if $\,(\pep\cup\pem)\cap\hf(\x)=$\hskip4pt\emp\hskip4.5ptfor its 
Gauss mapping $\,\hf$.

Relations (\ref{tre}) -- (\ref{tno}) remain valid in this case. In fact, the 
pull\-backs $\,\hf^*\la^+$ and $\,\hf^*\mu^+$ are canonically trivialized by 
the $\,\la^+$ and $\,\mu^+$ components of $\,\pi^*\psi$, so that 
(\ref{psg}\hs-b) yields (\ref{tre}), while (\ref{tnu}) now gives 
$\,\hf^*\tau=\hf^*\la^-=\hf^*\overline\nu$, and hence (\ref{tno}). (Note that 
$\,T\x=F^*\tau\,$ and $\,\nu_f=F^*\nu$, cf.\ end of \S\ref{libu}.)

Immersions of oriented real surfaces in spin\ce\mfd s with Higgs 
fields obviously include, as a special case, their immersions in \acsu s (see 
\S\ref{hggs}). However, for \tri s one can also turn this relation around and 
view the former as a special case of the latter. In fact, the image of such an 
immersion $\,f:\x\to\y\,$ must lie in the open submanifold 
$\,\y\nh\smallsetminus\zr\,$ on which the Higgs field is nonzero (cf.\ 
Remark~\ref{recpp}) and, by Theorem~\ref{thuno}, 
$\,f:\x\to\y\nh\smallsetminus\zr\,$ then is a \tri\ in the \acsu\ formed by 
$\,\y\nh\smallsetminus\zr\,$ with the residual \acst\ $\,J\,$ described in 
Remark~\ref{reres}.

Let $\,\psi\,$ again be a Higgs field for a \spst\ $\,(\sg^+\!,\sg^-\!,\ka)\,$ 
over a $\,4$\mfd\ $\,\y$, and let $\,E^+$ denote the total space of the 
unit-circle bundle of $\,\ka$. Every \tri\ $\,f:\x\to\y\,$ of an oriented 
closed real surface $\,\x\,$ gives rise to a mapping $\,\varXi:\x\to E^+$ with 
$\,\varXi(x)=(y,\rho/|\rho|)$, where, for any $\,x\in\x$, we set $\,y=f(x)\,$ 
and $\,\rho=\phi\hs\otimes\chi\,$ with $\,\phi,\chi\,$ standing for the 
$\,\la^+$ and $\,\mu^+$ components of $\,\pi^*\psi\,$ at $\,\hf(x)$. (Thus, 
$\,\phi\otimes\chi\in(\pi^*\ka)_{\hf(x)}=\ka_y$ by (\ref{psg}\hs-b).) The 
homotopy class of this mapping $\,\x\to E^+$ may be called the {\em oriented 
Maslov invariant\/} of $\,f$. (See~\cite{derdzinski-januszkiewicz} and 
Remark~\ref{rewed} below.) 

Here $\,\pi:\gy\to\y\,$ is, as usual, the bundle projection of the 
Grassmannian manifold, and $\,\hf\,$ is the Gauss mapping of $\,f$, 
cf.\ end of \S\ref{grop}.

According to Gromov~\cite[p.~192]{gromov}, the $\,h\hskip.3pt$-principle holds 
for \tri s of closed real surfaces in \acsu s. Therefore, the oriented Maslov 
invariant classifies such immersion up to the equivalence relation of being 
homotopic through \tri s. (See~\cite{derdzinski-januszkiewicz}.) 

The last conclusion remains valid, more generally, for \tri s of oriented 
closed real surfaces in a $\,4$\mfd\ endowed with a \spst\ and a Higgs-field. 
Namely, $\,\zr\,$ in Remark~\ref{reres} is a discrete subset of $\,\y$. Hence, 
for dimensional reasons, the set of homotopy classes of mappings $\,\x\to E^+$ 
remains unchanged when $\,E^+$ is replaced by its portion lying over 
$\,\y\nh\smallsetminus\zr$, while the set of equivalence classes of \tri s 
$\,\x\to\y\,$ is the same as for \tri s $\,\x\to\y\nh\smallsetminus\zr$.

\begin{rem}\label{rewed}We have $\,v\wedge w=-\hs2i\hs\phi\otimes\chi\,$ 
whenever $\,(v,w)\,$ is a pos\-i\-tive-o\-ri\-ented $\,g$-or\-tho\-nor\-mal 
basis of a plane $\,\mathcal{T}\subset T_y\y\,$ at a point $\,y\,$ of a 
$\,4$\mfd\ $\,\y\,$ with a Hermitian structure $\,(J,g)$, and 
$\,\phi,\hs\chi\,$ denote the $\,\la^+$ and $\,\mu^+$ components of 
$\,\pi^*\psi\,$ at $\,(y,\mathcal{T})\,$ relative to the decomposition 
(\ref{psg}\hs-a) for the \spst\ $\,(\sg^+\!,\sg^-\!,\ka)\,$ and the unit 
section $\,\psi\,$ of $\,\sg^+$ associated with $\,(J,g)\,$ as in 
\S\ref{acsu}. Here $\,v\wedge w\,$ is the {\em complex\/} exterior product, so 
that $\,v\wedge w\in\ka_y$ with $\,\ka=\detc\hs T\y\,$ (see \S\ref{acsu}); 
also, $\,\phi\otimes\chi\in(\pi^*\ka)_{(y,\mathcal{T})}=\ka_y$ in view of 
(\ref{psg}\hs-b).

Consequently, the above definition of the Maslov invariant is equivalent to 
that in~\cite{derdzinski-januszkiewicz}, which uses 
$\,v\wedge w/|v\wedge w|\,$ instead of our $\,\rho/|\rho|$. (The $\,-\hs i\,$ 
factor, representing a rotation, leaves the homotopy class unaffected.)

Relation $\,v\wedge w=2i\hs\phi\otimes\chi\,$ amounts to
$\,A\psi\wedge B\psi=-\hs2i\hs\phi\wedge\chi\,$ for any 
pos\-i\-tive-o\-ri\-ented orthonormal basis $\,(A,B)\,$ of 
$\,\mathcal{T}\in \gv$, where $\,\mathcal{V}\,$ is the space (\ref{vee}) for a 
fixed \spingeo\ $\,(\csg^+\!,\csg^-\!,\cka)$, and $\,\phi,\chi\,$ are the 
$\,\mathcal{L}\,$ and $\,\mathcal{L}^\perp$ components of $\,\psi\in\csg^+$, 
with complex lines $\,\mathcal{L}\in P(\csg^+)\,$ and 
$\,\mathcal{L}'\in P(\csg^-)\,$ such that 
$\,\mathcal{T}=\varPsi(\mathcal{L},\mathcal{L}')\,$ for $\,\varPsi\,$ as in 
(\ref{dif}). Note that it is the Clifford multiplication by $\,\psi(y)$, for 
$\,y\in\y$, that turns $\,\mathcal{V}=T_y\y\,$ into a {\em complex\/} vector 
space by identifying it with $\,\sg^-_y$ (see \S\ref{acsu}). 

Equality $\,A\psi\wedge B\psi=-\hs2i\hs\phi\wedge\chi\,$ holds since, due to 
the matrix form of any $\,A\in\mathcal{V}\,$ (see Remark~\ref{revee}), 
$\,A^*\!A\,$ equals $\,|A|^2$ times the identity operator of $\,\csg^+\nh$. 
The Euclidean inner product $\,\langle\,,\rangle\,$ of $\,\mathcal{V}\,$ thus 
is given by $\,4\langle A,B\rangle=\,\tracer\,\hs A^*\nh B\,$ which, if 
$\,A\ne0$, coincides with $\,|A|^{-2}\hs\tracer\,\hs A^{-1}\nh B$. However, 
$\,A^{-1}\nh B\,$ has the eigenvalues $\,z,\overline z\,$ for some 
$\,z\in\bbC\,$ with $\,\hs{\rm Im}\,z>0$, realized by eigenvectors in 
$\,\mathcal{L}\,$ and, respectively, $\,\mathcal{L}^\perp$. (See end of 
\S\ref{xpdi}.) As $\,\langle A,B\rangle=0$, that is, 
$\,\hs\tracer\,\hs A^{-1}\nh B=0$, we then have $\,\hs{\rm Re}\,z=0$. Also, 
$\,A,B\,$ are linear isometries (being homotheties of norm $\,1$, cf. 
Remark~\ref{revee}), and so $\,|z|=1$, i.e., $\,z=i$. Now $\,B\phi=iA\phi$, 
$\,B\chi=-\hs iA\chi$. Thus, from $\,\psi=\phi+\chi$, we get 
$\,i\hs A\psi\wedge B\psi=2A\phi\wedge A\chi=2\hs(\det A)\hs\phi\wedge\chi\,$ 
with $\,\det A=1\,$ (see Remark~\ref{revee}).
\end{rem}

\section{Pseu\-do\-hol\-o\-mor\-phic immersions}\label{pshi}
Let $\,\psi\,$ be a fixed Higgs field (cf.\ \S\ref{hggs}) for a \spst\ 
$\,(\sg^+\!,\sg^-\!,\ka)\,$ over a $\,4$\mfd\ $\,\y$. We say that an immersion 
$\,f:\x\to\y\,$ of an oriented real surface $\,\x\,$ is {\em 
pseu\-do\-hol\-o\-mor\-phic\/} if the image $\,\hf(\x)\subset\gy\,$ of the 
Gauss mapping of $\,f\,$ (see end of \S\ref{grop}) is contained in the set 
$\,\pep\,$ described in \S\ref{hggs}, that is, if every point of $\,\x\,$ is a 
$\,\pep$-com\-plex point of $\,f$.

In terms of the set $\,\zr\,$ of zeros of $\,\psi$, an immersion 
$\,f:\x\to\y\,$ of an oriented real surface is \psh\ if and only if its 
restriction to $\,\x\smallsetminus f^{-1}(\zr)\,$ is \psh, in the sense of 
\S\ref{cpts}, relative to the residual \acst\ $\,J\,$ on 
$\,\y'=\y\nh\smallsetminus\zr\,$ described in Remark~\ref{reres}. (In fact, by
Remark~\ref{recpp}, all points of $\,f^{-1}(\zr)\,$ are $\,\pep$-com\-plex 
points of $\,f$.) Such an immersion therefore induces an \acst\ on 
$\,\x\smallsetminus f^{-1}(\zr)\,$ which is always integrable as $\,\dim\x=2$, 
and has an extension to a complex structure on $\,\x\,$ (by Remark~\ref{reext} 
below, since $\,\zr\,$ is discrete).
\begin{rem}\label{reext}Let $\,g\,$ be a Riemannian metric on a manifold 
$\,\y\,$ and let $\,J\,$ be an \acst\ on $\,\y'=\y\smallsetminus\{y\}$, 
compatible with $\,g$, for some given point $\,y\in\y$. If $\,\x\subset\y\,$ 
is a two\diml\ real submanifold such that $\,y\in\x\,$ and the inclusion 
mapping $\,\x\hs'\to\y'$ of $\,\x\hs'=\x\smallsetminus\{y\}\,$ is 
\psh\  relative to $\,J$, then the (integrable) \acst\ 
which $\,J\,$ induces on $\,\x\hs'$ has a $\,C^\infty$ extension to $\,\x$.

This is clear as $\,J\,$ restricted to $\,T\x\hs'$ is uniquely determined, as 
in Remark~\ref{relin}, by the metric on $\,\x\hs'$ induced by $\,g\,$ and the 
orientation of $\,\x\hs'$ induced by $\,J\,$ via (\ref{ori}), both of which 
admit obvious extensions to $\,\x$.
\end{rem}

\section{A projection of $\,\cp^n$ onto $\,S^{\hs2n}$}\label{prcs}
Throughout this section $\,V\hs$ is a fixed complex vector space of complex 
dimension $\,n<\infty\,$ with a Hermitian inner product 
$\,\langle\,,\rangle\,$ and the corresponding norm $\,|\,\,|$, while $\,S,P\,$ 
stand for its natural compactifications: namely, $\,S\hs\approx\hs S^{2n}$ is 
obtained by adding to $\,V\hs$ a new point at infinity, that is, 
$\,S=V\cup\{\infty\}$, while $\,P\hs\approx\hs\cp^n$ is the set of all complex 
lines through zero in the direct sum (product) complex vector space 
$\,V\hskip-1.2pt\times\bbC\hs$, i.e., $\,P=P(V\hskip-1.2pt\times\bbC)\,$ (cf.\ 
\S\ref{prel}). We use the homogeneous-coordinate notation, so that 
$\,\lj y,z\rj\in P\,$ is the 
line spanned by $\,(y,z)\in V\hskip-1.2pt\times\bbC\,$ whenever $\,y\in V$, 
$\,z\in\bbC\,$ and $\,|y|+|z|>0$. Thus, $\,P=V\cup H$, where we treat $\,V\hs$ 
as a subset of $\,P$, identifying any $\,y\in V\,$ with $\,\lj y,1\rj$, and 
$\,H\hs\approx\hs\cp^{n-1}$ is the hyperplane of all $\,\lj y,0\hs\rj\,$ with 
$\,y\in V\smallsetminus\{0\}$. For $\,y\in V\smallsetminus\{0\}\,$ we now have 
$\,\lj y,0\hs\rj=L\times\{0\}\in H$, with $\,L=\bbC y\in P(V)$.

Finally, $\,\varPhi,\,\proj\hs\,$ denote the {\em inversion\/} 
$\,\varPhi:S\to S\,$ with $\,\varPhi(y)\,=\,y/|y|^2$ for 
$\,y\in V\smallsetminus\{0\}$, $\,\varPhi(0)=\infty$, $\,\varPhi(\infty)=0$, 
and the {\em projection\/} $\,\hs\proj:P\to S\,$ with 
$\,\hs\proj(\lj y,z\rj)=y/z\,$ if $\,z\ne0\,$ and 
$\,\hs\proj(\lj y,0\hs\rj)=\infty$.

The $\,C^\infty$-man\-i\-fold structure of $\,S\,$ is introduced by an atlas 
of two $\,V$-val\-ued charts, formed by 
$\,\hs{\rm Id}\hs:S\smallsetminus\{\infty\}=V\to V\,$ and 
$\,\varPhi:S\smallsetminus\{0\}\to V$. As one easily verifies in both charts, 
$\,\varPhi\,$ is a diffeomorphism $\,S\to S$, while $\,\hs\proj:P\to S\,$ is 
always of class $\,C^\infty$ and, if $\,n=1$, it is a diffeomorphism as well. 
(A function of $\,\lj y,z\rj\,$ is differentiable if and only if it is 
differentiable in $\,(y,z)$.) Thus, $\,\hs\proj:V\cup H\to V\cup\{\infty\}\,$ 
is the result of combining $\,\hs{\rm Id}\hs:V\to V\,$ with the constant 
mapping $\,H\to\{\infty\}$. In other words, $\,\hs\proj\hs\,$ identifies 
$\,S\,$ with the quotient of $\,P\,$ obtained by collapsing $\,H\,$ to the 
single point $\,\infty$.

For any $\,L\in P(V)$, i.e., any complex line $\,L\,$ through $\,0\,$ in $\,V$,
\begin{equation}
L\,{\rm\ is\ the\ image\ of\ the\ differential\ of\ }\,\varPhi\hs
\circ\hs\proj\hs\,{\rm\ at\ the\ point\ }\,L\times\{0\}\in P{\rm,}\label{iml}
\end{equation}
with $\,L\,$ treated as a real vector subspace of $\,V=T_0V=T_0S$. In fact, 
since $\,\hs\proj\hs\,$ is constant on $\,H$, its real rank at 
$\,L\times\{0\}\,$ is at most $\,2$, and so (\ref{iml}) will follow if we show 
that the image in question contains $\,L$. To this end, note that, as we just 
saw, $\,\hs\proj\hs\,$ and $\,\varPhi\,$ send the complex projective line 
$\,P\hs'=P(L\times\bbC)\subset P\,$ and the $\,2$-sphere 
$\,S\hs'=L\cup\{\infty\}\subset S\,$ diffeomorphically onto $\,S\hs'$. Thus, 
$\,\varPhi\hs\circ\hs\proj:P\hs'\to S\hs'$ is a diffeomorphism; hence, the 
image of its differential at $\,L\times\{0\}\,$ is $\,T_0S\hs'=T_0L=L$. Next, 
in terms of the decomposition $\,V=\bbC y\oplus y^\perp$ of $\,V=T_yV=T_yS\,$ 
for any $\,y\in V\smallsetminus\{0\}$,
\begin{equation}
{\rm the\ differential\ }\,d\varPhi_y{\rm\ is\ complex{\textrm-}linear\ on\ }
\,y^\perp{\rm\ and\ antilinear\ on\ }\,\bbC y\hs{\rm.}\label{dph}
\end{equation}
The $\,y^\perp$ part is obvious since $\,d\varPhi_y$ acts on the 
$\,\hs{\rm Re}\hskip1pt\langle\,,\rangle$-or\-thog\-o\-nal complement of 
$\,y$, which includes $\,y^\perp$, via multiplication by $\,|y|^{-2}$, as one 
sees applying $\,\varPhi\,$ to curves on which the norm $\,|\,\,|\,$ is 
constant. For the $\,\bbC y\,$ part, note that $\,L=\bbC y\,$ is tangent at 
$\,y\,$ to the $\,\varPhi$-invariant submanifold 
$\,L\smallsetminus\{0\}\subset V\hs$, which we may identify with 
$\,\bbC\smallsetminus\{0\}\,$ using the diffeomorphism of 
$\,\bbC\smallsetminus\{0\}\,$ onto $\,L\smallsetminus\{0\}\,$ sending $\,z\,$ 
to $\,zy/|y|$. This makes $\,\varPhi\,$ restricted to 
$\,V\smallsetminus\{0\}\,$ appear as the standard antiholomorphic inversion 
$\,z\mapsto z/|z|^2=1/\overline z$.

\section{Standard Higgs fields on the \hs4-sphere}\label{sphr}
We define a specific spin structure $\,(\sg^+\nh,\sg^-)\,$ over the sphere 
$\,S^{\hs4}$ (see end of \S\ref{spct}) by first identifying $\, S^{\hs4}$ with 
the quaternion projective line $\,\pw\approx\hp^1$ for a given quaternion 
plane $\,W\,$ (cf.\ \S\ref{prel} and Lemma~\ref{sfour}(a) below), and then 
choosing $\,\sg^+\nh,\sg^-$ to be, respectively, the tautological 
$\,\bbH$-line bundle over $\,\pw$, with the fibre $\,L\,$ at any $\,L\in \pw$, 
and the quotient $\,\sg^-=\theta/\sg^+$, where $\,\theta=\pw\times W\,$ is the 
product bundle over $\,\pw\,$ with the fibre $\,W$. A fixed quaternion inner 
product $\,\langle\,,\rangle\,$ in $\,W\,$ now turns both $\,\sg^\pm$ into 
normed quaternion line bundles, with the structure group $\,\hs\sug(2)\,$ 
(cf.\ the matrix in (\ref{mlt})), which also makes them Hermitian plane 
bundles (see the lines following (\ref{ori})), and yields the 
norm-pre\-serv\-ing identifications 
$\,\ka\,=\,[\sg^+]^{\wedge2}=\,[\sg^-]^{\wedge2}$, required in \S\ref{spct}, 
for $\,\ka=\pw\times\bbC\hs$. Finally, equality 
$\,T[\pw]=\,\hs\hh(\sg^+\nh,\sg^-)$, due to (\ref{tlp}), provides the Clifford 
multiplication by any $\,A\in T_L[\pw]$, $\,L\in \pw$. As $\,A:L\to W/L\,$ is 
$\,\bbH$-linear and $\,\dimh L=\dimh W/L=1$, the properties of $\,A\,$ named 
in (\ref{vee}) follow easily.

Every nonzero vector $\,u\in W\,$ now gives rise to a section $\,\psi\,$ of 
$\,\sg^+$ with $\,\psi(L)=\hs\proj_Lu\in L\,$ for $\,L\in \pw$, where 
$\,\hs\proj_L$ is the $\,\langle\,,\rangle$-or\-thog\-o\-nal projection onto 
$\,L$. Thus, $\,\psi(L)=\langle u,w\rangle\hs w/|w|^2$ if $\,L=\bbH w$. For 
any $\,u\in W\smallsetminus\{0\}\,$ this $\,\psi\,$ is a Higgs field for 
the spin structure $\,(\sg^+\nh,\sg^-)\,$ over $\,\pw$, cf.\ \S\ref{hggs}, and 
will be called a {\em standard Higgs field\/} for $\,(\sg^+\nh,\sg^-)$. In 
fact, $\,\psi\,$ has just one zero, at $\,u^\perp\in \pw$, and is 
transverse to the zero section by (i) in \S\ref{trsv}.
\begin{lem}\label{sfour}Let\/ $\,\psi\,$ be the standard Higgs field for the 
spin structure\/ $\,(\sg^+\nh,\sg^-)\,$ over\/ $\,\pw\,$ corresponding to a 
nonzero vector\/ $\,u\,$ in a quaternion inner-product plane\/ $\,W$, and let 
the quaternion line\/ $\,V=u^\perp\in \pw\,$ be treated as a complex Hermitian 
plane, cf.\ {\rm\S\ref{prel}}. If\/ $\,S=V\cup\{\infty\}\,$ is the\/ 
$\,4$-sphere described in\/ {\rm\S\ref{prcs}} and\/ $\,\vt:S\to \pw\,$ is the 
extension, with\/ $\,\hs\vt(\infty)=V\nh$, of\/ $\,\vt:V\to \pw\,$ defined 
in\/ {\rm Remark~\ref{redvt}}, then
\begin{enumerate}
  \def\theenumi{{\rm(\alph{enumi}}}
\item $\vt:S\to \pw\,$ is a\/ $\,C^\infty$ diffeomorphism.
\item The residual \acst\/ $\,J\,$ on\/ $\,\pw\smallsetminus\{V\}$, 
associated with\/ $\,\psi\,$ as in\/ {\rm Remark~\ref{reres}}, corresponds 
under\/ $\,\vt\,$ to the obvious structure on the complex vector space\/ 
$\,V=S\smallsetminus\{\infty\}$.
\end{enumerate}
\end{lem}
\begin{proof}Both $\,\vt\,$ and $\,\vt^{-1}$ appear as 
real-rational mappings in standard projective coordinates for $\,\pw\,$ 
(cf.\ Remark~\ref{redvt}) and the two charts for $\,S\,$ (see \S\ref{prcs}), 
which proves (a). Given $\,y\in V\hs$ and $\,v\in V=T_yV=T_yS$, the operator 
$\,A:L\to W/L\,$ with $\,L=\bbH w\,$ for $\,w=y+u$, corresponding to 
$\,\dvt_yv\,$ under (\ref{tlp}), sends $\,w=y+u\,$ to the coset $\,v+L\,$ (see 
Remark~\ref{redvt}), so that, by $\,\bbH$-linearity of $\,A$, the $\,A$-image 
of $\,\psi(L)=\hs\proj_Lu=\langle u,w\rangle\hs w/|w|^2$ is the coset of 
$\,L\,$ in $\,W\,$ containing $\,\langle u,w\rangle\hs v/|w|^2$, i.e., $\,v\,$ 
times the {\em real\/} scalar $\,|u|^2/(|y|^2+|u|^2)$. (Note that 
$\,\langle u,w\rangle=|u|^2$ as $\,y\in V=u^\perp$.) The dependence of that 
coset on $\,v\,$ thus is $\,\bbC$-linear, relative to the obvious complex 
structure of $\,V$. This yields (b), completing the proof.
\end{proof}
\begin{rem}\label{repsh}Combined with the second paragraph of \S\ref{pshi}, 
Lemma~\ref{sfour} shows that an immersion $\,f\,$ of an oriented real surface 
$\,\x\,$ in the $\,4$-sphere $\,S=V\cup\{\infty\}\,$ is \psh\ relative to some 
standard Higgs field in the spin structure $\,(\sg^+\nh,\sg^-)\,$ over $\,S\,$ 
if and only if, restricted to $\,\x\smallsetminus f^{-1}(\infty)$, it is \psh\ 
as an immersion into the complex plane $\,V$.
\end{rem}

\section{Smoothness at infinity for holomorphic curves in $\,\bbC^{\hs n}$}
\label{smth}
As a step toward a classification proof in \S\ref{rsps}, we will now show that 
any smooth real surface in the sphere 
$\,S^{\hs2n}=\bbC^{\hs n}\cup\{\infty\}$, obtained by adding the point 
$\,\infty\,$ to a holomorphic curve in $\,\bbC^{\hs n}\nh$, is the image of a 
holomorphic curve in $\,\cp^n=\bbC^{\hs n}\cup\cp^{n-1}$ under a natural 
projection (see \S\ref{prcs}) which restricted to $\,\bbC^{\hs n}$ is the 
identity, and sends the hyperplane $\,\cp^{n-1}$ onto $\,\infty$. 
\begin{lem}\label{pshol}Given\/ 
$\,V,\langle\,,\rangle,S,\infty,P,H,\varPhi,\,\proj\hs\,$ as in\/ 
{\rm\S\ref{prcs}}, let\/ $\,\x\,$ be a two\diml\ real submanifold of\/ $\,S\,$ 
such that\/ $\,\infty\in\x\,$ and\/ $\,\x\smallsetminus\{\infty\}\,$ is a 
complex submanifold of\/ $\,V=S\smallsetminus\{\infty\}$. Furthermore, let\/ 
$\,L=d\varPhi_\infty(T_\infty\x)$ be the real vector subspace of 
$\,V=T_0V=T_0S\hs$ obtained as the image of\/ $\,T_\infty\x\,$ under the 
differential of\/ $\,\varPhi\,$ at\/ $\,\infty$. Then
\begin{enumerate}
  \def\theenumi{{\rm(\roman{enumi}}}
\item $L\,$ is a complex line in\/ $\,V$.
\item $\x\,$ is the\/ $\,\hs\proj$-image of a one\diml\ complex submanifold\/ 
$\,\hat\x\,$ of\/ $\,P\,$ which intersects\/ $\,H$, transversally, at the 
single point\/ $\,L\times\{0\}$.
\end{enumerate}
\end{lem}
A proof of Lemma~\ref{pshol} is given in the next section.

\begin{rem}\label{reima}Given a complex inner-product space $\,V\hs$ of any 
dimension $\,n<\infty$, we can now describe all two\diml\ real 
submanifolds $\,\x\,$ of the $\,2n$-sphere $\,S=V\cup\{\infty\}\,$ for which 
$\,\infty\in\x\,$ and $\,\x\smallsetminus\{\infty\}\,$ is a complex 
submanifold of $\,V$. Namely, they coincide with the $\,\hs\proj$-images of 
the one\diml\ complex submanifolds $\,\hat\x\,$ of the projective space 
$\,P=V\cup H\,$ that intersect the hyperplane $\,H\hs\approx\hs\cp^{n-1}\nh$, 
transversally, at a single point. In fact, every such $\,\x\,$ is of this form 
by Lemma~\ref{pshol}(ii), while the converse statement is immediate since 
$\,\hs\proj:\hat\x\to S\,$ is an immersion: by (\ref{iml}), its real rank at 
the intersection point of $\,\hat\x\,$ with $\,H\,$ equals $\,2$.
\end{rem}

\section{Proof of Lemma~\ref{pshol}}\label{proo}
Denoting $\,\vg\,$ the unit circle about $\,0\,$ in $\,L$, let us choose 
$\,\ve\in(0,\infty)\,$ and a $\,C^\infty$ mapping 
$\,[\hs0,\ve]\times \vg\ni(t,u)\mapsto y(t,u)\in N\nh\cap V$, where 
$\,N=\varPhi(\x)\subset S=V\cup\{\infty\}$, such that $\,y(0,u)=0$, 
$\,\dot y(0,u)=u\,$ with $\,\dot y(t,u)=\hs d\hs[y(t,u)]/dt$, and 
$\,y(t,u)\ne0\,$ unless $\,t=0$, for every $\,u\in \vg$. For instance, we 
might set $\,y(t,u)=\exp_{\hs0}tu$, using the exponential mapping at 
$\,0\in N\,$ of any Riemannian metric on $\,N$.

To prove (i), let us fix $\,u\in \vg\,$ and, for any given $\,t\in(0,\ve]$, 
write $\,x=y(t,u)$, $\,v=\dot y(t,u)$. Thus, 
$\,d\varPhi_xv\in T_{\varPhi(x)}\x$, as $\,v\in T_xN\,$ and $\,N=\varPhi(\x)$. 
Since $\,T_{\varPhi(x)}\x\,$ is closed under multiplication by $\,-\hs i$, 
combining this with (\ref{dph}) we obtain $\,i\vrad-i\vtng\in T_xN$, where 
$\,\vrad,\vtng$ are the components of $\,v\,$ relative to the decomposition 
$\,V=\bbC x\oplus x^\perp$. In the last relation, $\,\vrad,\vtng$ and $\,x\,$ 
all depend on $\,t$, and taking its limit as $\,t\to0\,$ we get 
$\,iu\in T_0N$, since $\,\vrad\to u$, $\,\vtng\to0\,$ as $\,t\to0\,$ (see the 
next paragraph) and, clearly, $\,x\to x(0)=0$. Therefore, $\,u,iu\in L=T_0N$, 
which yields (i).

To obtain relations $\,\vrad\to u\,$ and $\,\vtng\to0\,$ as $\,t\to0$, note 
that $\,v\to u\,$ by continuity of $\,\dot y(t,u)\,$ in $\,(t,u)$. (In all 
limits, $\,t\to0$.) Next, $\,x/t\to u\,$ (since $\,\dot y(0,u)=u$), and so 
$\,x/|x|\to u$, as $\,x/|x|=(x/t)/|x/t|$. Thus, 
$\,\vrad=\langle v,x/|x|\rangle\hs x/|x|\to\langle u,u\rangle\hs u=u$, and 
$\,\vtng=v-\vrad\to u-u=0$.

Since $\,\lj y,z\rj=[t^{-1}y,t^{-1}z\rj\,$ for $\,\lj y,z\rj\in P\,$ and 
$\,t\in\bbR\smallsetminus\{0\}$, the definitions of $\,\varPhi\,$ and 
$\,\hs\proj\hs\,$ give 
$\,\hs\proj^{-1}(\varPhi(y(t,u)))=\lj t^{-1}y(t,u),t^{-1}|y(t,u)|^2\rj$, while 
$\,t^{-1}y(t,u)\to u\,$ and $\,t^{-1}|y(t,u)|^2\to0\,$ as $\,t\to0$, uniformly 
in $\,u\in \vg$.

In fact, from Taylor's formula, $\,y(t,u)=tu+t^2\hs\omega(t,u)\,$ with a 
$\,C^\infty$ mapping $\,\hs\omega:[\hs0,\ve]\times\vg\to V$. (Integration by 
parts gives $\,\hs\omega(t,u)=\int_0^1(1-s)\hs\ddot y(st,u)\,ds$.) As 
$\,y(t,u)=t(u+t\hs\omega(t,u))$, convergence follows from continuity of 
$\,\hs\omega$, and is uniform due to compactness of $\,\vg$.

In other words, the restriction to 
$\,\x\smallsetminus\{\infty\}=\varPhi(N\smallsetminus\{0\})\,$ of the 
biholomorphism 
$\,\hs\proj^{-1}:S\smallsetminus\{\infty\}\to P\smallsetminus H\,$ (i.e., of 
the identity mapping of $\,V$) has a limit at $\,\infty\,$ equal to 
$\,L\times\{0\}\in H$, that is, to $\,\lj u,0\hs\rj\,$ for any 
$\,u\in \vg$. Setting $\,\hat f=\hs\proj^{-1}$ on 
$\,\x\smallsetminus\{\infty\}\,$ and $\,\hat f(\infty)=L\times\{0\}\,$ we now 
obtain a continuous mapping $\,\hat f:\x\to P$, which is holomorphic on 
$\,\x\smallsetminus\{\infty\}$. Also, the complex structure of 
$\,\x\smallsetminus\{\infty\}\,$ has an extension to $\,\x$.

Such an extension exists in view of Remark~\ref{reext} applied to 
$\,\y=S\smallsetminus\{0\}$, $\,y=\infty$, and $\,J\,$ which is the standard 
integrable \acst\ on $\,M'=V\smallsetminus\{0\}=S\smallsetminus\{0,\infty\}$, 
along with the metric $\,g\,$ obtained as the pull\-back under $\,\varPhi\,$ 
of the standard Euclidean metric $\,\hs{\rm Re}\hskip1pt\langle\,,\rangle\,$ 
on $\,V$. Note that, since $\,\varPhi\,$ is conformal, $\,g\,$ equals a 
positive function times $\,\hs{\rm Re}\hskip1pt\langle\,,\rangle$, and so is 
still compatible with $\,J$.

Some neighborhood of $\,\infty\,$ in $\,\x\,$ is therefore biholomorphic 
to a disk in $\,\bbC\hs$, so that the mapping $\,\hat f:\x\to P$, continuous 
on $\,\x\,$ and holomorphic on $\,\x\smallsetminus\{\infty\}$, must be 
holomorphic on $\,\x$. Also, by continuity, $\,\hs\proj\circ\hat f\,$ 
coincides with the inclusion mapping $\,\x\to S$. Thus, $\,\hat f\,$ is a 
$\,C^\infty$ embedding, and hence a holomorphic embedding. Setting 
$\,\hat\x=\hat f(\x)\,$ we thus obtain a biholomorphism 
$\,\hat f:\x\to\hat\x$, the inverse of which is $\,\hs\proj\hs\,$ restricted 
to $\,\hat\x\,$ (and valued in $\,\x$). Since $\,\hs\proj\hs\,$ is constant on 
$\,H$, transversality of $\,\hat\x\,$ and $\,H\,$ follows, completing the 
proof of Lemma~\ref{pshol}.

\section{Real surfaces, pseu\-do\-hol\-o\-mor\-phic\-al\-ly immersed in 
$\,S^{\hs4}$}\label{rsps}
For the $\,4$-sphere $\,S=V\cup\{\infty\}\,$ obtained as in \S\ref{prcs} (with 
$\,n=2$) from a Hermitian plane $\,V,$ an immersion $\,f:\x\to S\,$ of an 
oriented real surface $\,\x\,$ will be called {\em \psh\/} if its restriction 
to $\,\x\smallsetminus f^{-1}(\infty)\,$ is a \psh\ immersion in the complex 
plane $\,V$. According to Remark~\ref{repsh}, this amounts to its being \psh\ 
for a certain standard Higgs field in the spin structure 
$\,(\sg^+\nh,\sg^-)\,$ over $\,S\,$ defined as in \S\ref{sphr}.

In this section we classify such immersions of {\em closed\/} surfaces $\,\x$, 
beginning with the special case of embeddings.
\begin{thm}\label{thdue}Let\/ $\,P\,\approx\,\cp^2$, $\,S\,\approx\,S^{\hs4}$ 
and\/ $\,\hs\proj:P\to S\,$ be defined as in\/ {\rm\S\ref{prcs}} with\/ 
$\,n=2$, for some two\diml\ complex inner-product space\/ $\,V$, so that\/ 
$\,P=V\cup H\,$ and\/ $\,S=V\cup\{\infty\}$, where\/ $\,H\,\approx\,\cp^1$ is 
a projective line in\/ $\,P$, while $\,\hs\proj\hs\,$ sends $\,H\,$ to 
$\,\infty\,$ and\/ coincides with the identity on\/ $\,V$.

The oriented closed real surfaces pseu\-do\-hol\-o\-mor\-phic\-al\-ly embedded 
in\/ $\,S$ then are all diffeomorphic to\/ $\,S^2$, and coincide with the\/ 
$\,\hs\proj$-images of complex projective lines other than\/ $\,H\,$ in the 
projective space\/ $\,P$.
\end{thm}
In fact, $\,\infty\in\x\,$ for any pseu\-do\-hol\-o\-mor\-phic\-al\-ly 
embedded, oriented closed real surface $\,\x\subset S\,$ (as $\,\x\,$ cannot 
lie entirely in $\,V=S\smallsetminus\{\infty\}$), and so Remark~\ref{reima} 
with $\,n=2\,$ gives $\,\x=\hs\proj(\hat\x)\,$ for a one\diml\ compact complex 
submanifold $\,\hat\x\,$ which intersects the projective line 
$\,H$, transversally, at a single point. By Chow's theorem 
\cite{griffiths-harris}, $\,\hat\x\,$ is algebraic, and hence of degree one, 
as required. The converse statement is clear from Remark~\ref{reima}.
\begin{thm}\label{thtre}Under the hypotheses of\/ {\rm Theorem~\ref{thdue}}, 
the \psh\ immersions\/ $\,f:\x\to S\,$ of any oriented real surface\/ $\,\x\,$ 
are nothing else than the composites\/ $\,\hs\proj\hs\circ\hat f$, where\/ 
$\,\hat f\,$ runs through all \psh\ immersions\/ $\,\x\to P\,$ transverse to\/ 
$\,H$.
\end{thm}
This is immediate if one replaces $\,\x\,$ by the $\,f$-image of a suitable 
neighborhood in $\,\x\,$ of any given point of $\,f^{-1}(\infty)\,$ and, 
again, applies Remark~\ref{reima}.

Unlike Theorem~\ref{thdue}, the last result does not even assume closedness of 
$\,\x$. However, when the surface $\,\x\,$ {\em is\/} closed, Chow's theorem 
\cite[p.~167]{griffiths-harris} implies that the image $\,\hat f(\x)\,$ is an 
algebraic curve of some degree $\,d$. In the simplest case, where some finite 
number $\,\delta\,$ of ordinary double points constitute both the only 
singularities of $\,\hat f(\x)\,$ and the only self-intersections of $\,f$, 
the genus of $\,\x\,$ is $\,[(d-1)(d-2)-2\delta\hs]/2$, by Pl\"ucker's formula 
\cite[p.~280]{griffiths-harris}.

\section{Blow-ups and blow-downs involving Higgs fields}\label{blow}
The one-to-one correspondence, in Theorem~\ref{thtre}, between \psh\ 
immersions $\,f:\x\to S\,$ on the one hand, and $\,\hat f:\x\to P\,$ on the 
other, is a special case of two mutually inverse constructions. They involve 
an embedded $\,\cp^1$ with the self-intersection number $\,+\hs1$, not 
$\,\hs-1$, and hence should not be confused with their obvious analogues in 
the complex category.

The first construction is {\em blow-down}. It may be applied to a complex 
surface $\,\hat\y\,$ along with a fixed complex submanifold 
$\,H\subset\hat\y$, biholomorphic to $\,\cp^1$, such that some neighborhood 
$\,\,\hat U\,$ of $\,H\,$ in $\,\hat\y\,$ has a biholomorphic identification 
with a neighborhood of a projective line $\,\cp^1$ in $\,\cp^2$, under which 
$\,H=\cp^1$. Contracting $\,H=\cp^1$ to a single point, denoted $\,\infty$, we 
now transform $\,\,\hat U\,$ into a neighborhood $\,\,U\,$ of $\,\infty\,$ in 
the four-sphere $\,S^{\hs4}=\bbC^2\cup\{\infty\}$. (This is the $\,n=2\,$ case 
of the construction in \S\ref{prcs}.) The replacement of $\,H\,$ by 
$\,\infty\,$ in $\,\hat\y\,$ thus leads to a new $\,4$\mfd\ $\,\y\,$ along 
with an $\,H$-collapsing projection $\,\hs\proj:\hat\y\to\y$. Smoothness of 
$\,M\,$ follows since $\,M\,$ is the result of gluing 
$\,\hat\y\smallsetminus H\,$ and $\,\,U\,$ together using the diffeomorphism 
$\,\hs\proj:\hat U\smallsetminus H\to\,U\smallsetminus\{\infty\}\,$ of 
\S\ref{prcs}. In addition, any Hermitian metric $\,g\,$ on $\,\hat\y\,$ can 
obviously be modified so as to coincide, near $\,H$, with the standard 
Fubini-Study metric on $\,\cp^2$. Our gluing procedure then may clearly be 
extended to the \spst s with Higgs fields over 
$\,\,\hat U\smallsetminus H\subset\cp^2$ and 
$\,\,U\smallsetminus\{\infty\}\subset S^{\hs4}$, of which the former 
corresponds to the (almost) Hermitian structure, cf.\ \S\ref{hggs}, and the 
latter involves a standard Higgs field on $\,S^{\hs4}$, with a zero at 
$\,\infty\,$ (see \S\ref{sphr}). The result is a \spst\ over $\,\y\,$ with a 
Higgs field vanishing at $\,\infty$.

The opposite construction is {\em blow-up}, for which the starting point is 
a \spst\ over a $\,4$\mfd\ $\,\y\,$ with a Higgs field that vanishes at some 
given point of $\,\y$, denoted $\,\infty$. An additional assumption is that 
on some neighborhood $\,\,U\,$ of $\,\infty\,$ in $\,\y\,$ the \spst\ and the 
Higgs field may be identified, via a suitable bundle isomorphism, with 
a standard Higgs field on $\,S^{\hs4}$, restricted to an open set (also 
denoted $\,\,U$) and having a zero at a point $\,\infty\in\hs U$. As in 
\S\ref{prcs}, $\,\,U\,$ is obtained from an open set $\,\,\hat U\subset\cp^2$ 
by collapsing a projective line $\,H=\cp^1\subset\,\hat U\,$ to the point 
$\,\infty$. Let $\,\hat\y\,$ now be the set obtained from $\,\y\,$ by 
replacing $\,\{\infty\}$ with $\,H$. For reasons similar to those in the 
preceding paragraph, $\,\hat\y\,$ is a smooth $\,4$\mfd\ with a $\,C^\infty$ 
projection mapping $\,\hs\proj:\hat\y\to\y\,$ sending $\,H\,$ onto 
$\,\{\infty\}\,$ and diffeomorphic (equal to the identity) on 
$\,\hat\y\smallsetminus H=\y\smallsetminus\{\infty\}$. Also, the original 
\spst\ and Higgs field have extensions from 
$\,\y\smallsetminus\{\infty\}=\hat\y\smallsetminus H\,$ to $\,\hat\y\,$ such 
that, on $\,\,\hat U$, they are the ones associated with the complex structure 
and the Fubini-Study metric of $\,\cp^2$.

Note that both constructions can be performed repeatedly, as long as the 
choices of $\,H\,$ (or, $\,\infty$) that we use are pairwise disjoint (or, 
distinct). In addition, one may also apply the standard {\em complex\/} 
blow-down or blow-up to a copy of $\,\cp^1$ biholomorphically embedded in 
$\,\hat\y\,$ with the self-intersection number $\,-\hs1\,$ or, respectively, 
to a point in $\,\y\,$ at which the Higgs field is nonzero and in a 
neighborhood of which the residual almost complex structure (see 
Remark~\ref{reres}) is integrable.

Finally, in both cases described above, for any oriented real surface $\,\x$, 
relation $\,f=\hs\proj\hs\circ\hat f\,$ defines a bijective correspondence 
between \psh\ immersions $\,\hat f:\x\to\hat\y\,$ which are 
transverse to $\,H$, and \psh\ immersions $\,f:\x\to\y$. This conclusion, 
obtained from Remark~\ref{reima} exactly as in our two-line ``proof'' of 
Theorem~\ref{thtre}, includes Theorem~\ref{thtre} as a special case with 
$\,\y=S\,$ and $\,\hat\y=P$.

\end{document}